\documentclass[12pt,draftcls,onecolumn]{IEEEtran}

\usepackage{amsfonts}
\usepackage{amssymb}
\usepackage{amsmath}
\usepackage{mathptmx}
\usepackage{cancel}
\usepackage{pstricks}
\usepackage{psfrag}
\usepackage{mathrsfs}
\usepackage{graphicx}
\usepackage{pgf,tikz}
\usepackage{algorithm}
\usepackage{algpseudocode}
\usepackage{soul}
\def\be{\begin{equation}}
\def\ee{\end{equation}}
\def\bea{\begin{eqnarray}}
\def\eea{\end{eqnarray}}
\def\beann{\begin{eqnarray*}}
\def\eeann{\end{eqnarray*}}

\newcommand{\rank}{{\rm rank}}

\def\ns{\hspace{-1mm}}
\newcommand{\real}{{\mathbb{R}}}

\def\spacingset#1{\def\baselinestretch{#1}\small\normalsize}
\spacingset{1.35}

\newtheorem{lemma}{Lemma}
\newtheorem{theorem}{Theorem}
\newtheorem{remark}{Remark}
\newtheorem{corollary}{Corollary}
\newtheorem{proposition}{Proposition}

\newtheorem{problem}{Problem}
\newtheorem{example}{Example}[section]
\newtheorem{assumption}{Assumption}[section]

\def\be{\begin{equation}}
\def\ee{\end{equation}}
\def\bea{\begin{eqnarray}}
\def\eea{\end{eqnarray}}
\def\beann{\begin{eqnarray*}}
\def\eeann{\end{eqnarray*}}

\def\ns{\hspace{-1mm}}

\def\proof{\noindent{\bf{\em Proof:}\ \ }}
\def\QED{\mbox{\rule[0pt]{1.5ex}{1.5ex}}}
\def\endproof{\hspace*{\fill}~\QED\par\endtrivlist\unskip}
\def\endex{\hspace*{\fill}~\Box\par\endtrivlist\unskip}

\newcommand{\ima}{\operatorname{im}}
\newcommand{\normrank}{\operatorname{normrank}}
\newcommand{\diag}{\operatorname{diag}}
\newcommand{\defi}{\stackrel{\text{\tiny def}}{=}}

\newcommand{\complex}{{\mathbb{C}}}

\def\gB{{\cal B}}

\def\gD{{\cal D}}

\def\gJ{{\cal J}}

\def\gL{{\cal L}}

\def\gP{{\cal P}}

\def\gR{{\cal R}}
\def\gS{{\cal S}}

\def\gU{{\cal U}}
\def\gV{{\cal V}}

\def\gX{{\cal X}}
\def\gY{{\cal Y}}
\def\gZ{{\cal Z}}

\def\bmat{\left[ \begin{array}}
\def\emat{\end{array} \right]}
\def\endex{\hspace*{\fill}~$\square$\par\endtrivlist\unskip}
\def\bmat{\left[ \begin{array}}
\def\emat{\end{array} \right]}
\def\bsmat{\left[ \begin{smallmatrix}}
\def\esmat{\end{smallmatrix} \right]}

\def\l{{\lambda}}
\def\gP{{\cal P}}
\def\gB{{\cal B}}
\def\gU{{\cal U}}
\def\gL{{\cal L}}
\def\gR{{\cal R}}

\def\gV{{\cal V}}

\def\gS{{\cal S}}

\def\gX{{\cal X}}
\def\i{{i}}
\newcommand{\spanR}{\operatorname{span}}

\begin{document}
\title{\LARGE{A structural solution to the \\ monotonic tracking control problem}}

\author{Lorenzo~Ntogramatzidis, Jean-Fran\c{c}ois Tr{\'e}gou{\"e}t, \\ Robert Schmid and Augusto Ferrante
 \thanks{L. Ntogramatzidis is with the Department of Mathematics and
Statistics, Curtin University, Perth,
Australia. E-mail: {\tt L.Ntogramatzidis@curtin.edu.au. }}

 \thanks{Jean-Fran\c{c}ois Tr{\'e}gou{\"e}t is with Universit\'{e} de Lyon, Laboratoire Amp\`{e}re CNRS UMR 5005, INSA-Lyon; F-69621, Villeurbanne, France. E-mail: {\tt jean-francois.tregouet@insa-lyon.fr. } (Research partially carried out at Curtin University).}

\thanks{R. Schmid is with the Department of Electrical and Electronic Engineering,
The University of Melbourne, Parkville, VIC 3010, Australia. E-mail: {\tt rschmid@unimelb.edu.au. }}

 \thanks{A. Ferrante is with the Dipartimento di Ingegneria dell'Informazione, Universit\`a di Padova, 
via Gradenigo, 6/B -- I-35131 Padova, Italy. E-mail: {\tt augusto@dei.unipd.it. }}

\thanks{Partially supported by the Australian Research Council under the grant FT120100604.}
}

\markboth{DRAFT}{Shell \MakeLowercase{\textit{et al.}}: Bare Demo of IEEEtran.cls for Journals}

\maketitle

\vspace{-1cm}

\IEEEpeerreviewmaketitle

\begin{abstract}
In this paper we present a method for designing a linear time invariant (LTI) state-feedback controller to monotonically track a constant step reference at any desired rate of convergence for any initial condition.
Necessary and sufficient constructive conditions are given to deliver a monotonic step response from all initial conditions. This method is developed for multi-input multi-output (MIMO) systems, and can be applied to square and non-square systems, strictly proper and non-strictly proper systems, and, importantly, also minimum and non-minimum phase systems. The framework proposed here shows that for MIMO LTI systems the objectives of achieving a rapid settling time, while at the same time
avoiding overshoot and/or undershoot, are not necessarily competing objectives.
\end{abstract}


\section{Introduction}
\label{secintro}
The problem of improving the shape of the step response curve for linear time invariant (LTI) systems is as old as control theory. Its relevance is seen in countless applications such as heating/cooling systems, elevator and satellite positioning, automobile cruise control and the positioning of a CD disk read/write head. The common element in these problems involves designing a control input for the system to make the output take a certain desired target value, and then keep it there.
 
 A fundamental issue in classical feedback control is the design of control laws that provide good performance both at steady state and during the transient. The \textit{steady state performance} is typically assumed to be satisfactory if, once the transient vanishes, the output of the system is constant and equal (or very close) to the desired value.
When dealing with the \textit{transient performance}, one is usually concerned with the task of reducing both the overshoot and the undershoot, or, ideally, of achieving a monotonic response that rapidly converges to the steady-state value.
 It is commonly understood that the objectives of achieving a rapid (short) settling time, while at the same time avoiding overshoot and undershoot, are competing objectives in the controller design, and must be dealt with by seeking a trade-off, see e.g. \cite{Franklin-PE-94,Dorf-B-08}, or any standard textbook on the topic.
While this is certainly the case for single-input single-output (SISO) systems, the control methods we develop and implement in this paper challenge this widely-held perception for the multi-input multi-output (MIMO) case. We show in particular that in the case of LTI MIMO systems, it is possible to achieve arbitrarily fast settling time and also a monotonic step response in all output components for any initial condition, which naturally implies the avoidance of overshoot/undershoot even in the presence of non-minimum phase invariant zeros. 

In contrast with the extensive literature for SISO systems, which includes -- but is far from being limited to -- \cite{Hoagg-B-07,Lin-F-97,Stewart-D-06,Darbha-B-03,Bement-J-04,Darbha-03,Darbha-B-02,Middleton-91,Lau-MB-03} and the references cited therein, to date there have been very few papers offering analysis or design methods for avoiding undershoot or overshoot in the step response of MIMO systems, see e.g. \cite{Johansson-02} and the references therein. 
The most famous among the classical methods that deal with general tracking control problems is the so-called {\em model matching problem}, see e.g. \cite{Wolovich-72, Malabre-K-84,Kucera-91}, which does not, in general, yield solvability conditions expressed solely in terms of the system's structure.

A recent contribution offering design methods for MIMO systems is \cite{Schmid-N-09}, where a procedure is proposed for the design of a state-feedback controller to yield a non-overshooting step response for LTI MIMO systems. Importantly, this design method is applicable to non-minimum phase systems, does not assume that the system state is initially at rest, and can be applied to both continuous-time and discrete-time (and also strictly proper or non-strictly proper) systems. Very recently it has been shown in \cite{Schmid-N-09-1} how the method can be adapted to obtain a non-undershooting step response. The key idea behind the approach in \cite{Schmid-N-09} and \cite{Schmid-N-09-1} is to design the feedback matrix that achieves the desired closed-loop eigenstructure in such a way that only a small number of the closed-loop system modes appear in each component of the tracking error (which is defined as the difference between the system output and the desired target value). Indeed, if the closed-loop eigenstructure can be constrained in such a way that each component of the tracking error is driven only by a single real-valued closed-loop mode -- which is an exponential in the form $e^{\lambda\,t}$ in the continuous time or a power term $\lambda^k$ in the discrete time -- the output of the system is monotonic in each output component regardless of the initial condition of the system, and hence both overshoot and undershoot are avoided. For systems where the closed-loop eigenstructure can be constrained so that the error involves only the sum of two or three exponential terms (or powers in the discrete case) in each component, the design method of \cite{Schmid-N-09} offers a search algorithm for the selection of suitable closed-loop modes that ensures that the step response is non-overshooting, non-undershooting, or monotonic from any given initial condition and target reference. 

A key limitation of the design methods given in \cite{Schmid-N-09} and \cite{Schmid-N-09-1} is the lack of analytic conditions, expressed in terms of the system structure, that guarantee the existence of a state-feedback controller that can deliver the desired transient response. In other words, the method of \cite{Schmid-N-09} and \cite{Schmid-N-09-1} does not provide a structural criterion to decide if the problem is solvable in terms of the problem data, nor does it guarantee that when the aforementioned matrix is singular one is allowed to conclude that the problem of achieving a monotonic response from any initial condition cannot be solved.
Moreover, as aforementioned, the design method involves a search for suitable closed-loop eigenvalues, and while this search can be conducted efficiently, the authors were unable to give any conditions guaranteeing a satisfactory search outcome. 
 The objective of this paper is to completely revisit the design method of \cite{Schmid-N-09} and \cite{Schmid-N-09-1} to the end of developing
 conditions expressed in terms of the 
  system structure that are necessary and sufficient to guarantee that the design method will deliver a state-feedback controller that yields a monotonic step response from any initial condition and for any constant reference signal. When this goal is achievable, we say that the control yields a {\em globally monotonic} response, by which we mean that the same feedback matrix yields a monotonic response from all initial conditions, and with respect to all possible step references. 
 
Thus, in this paper, for the first time in the literature, a complete and exhaustive answer to the problem of achieving a globally monotonic step response for a MIMO LTI system is provided. We show that for MIMO LTI systems the presence of non-minimum phase invariant zeros does not prevent a globally monotonic step response to be achievable. Indeed, even in the presence of one or more non-minimum phase invariant zeros, if the feedthrough matrix is allowed to be non-zero, it may still be possible to achieve a monotonic step response from any initial condition and for any constant reference signal. 

In the last part of the paper, we also offer a complete parameterisation of all the feedback matrices that achieve global monotonicity, thus opening the door to the formulation of optimisation problems whose goal is to exploit the available freedom to address further objectives such as minimum gain or improved robustness of the closed-loop eigenstructure in the same spirit of \cite{SICON}.\\[-3mm]

{\bf Notation.} In this paper, the symbol $\{0\}$ stands for the origin of a vector space. For convenience, a linear mapping $A: \gX \longrightarrow \gY$ between finite-dimensional vector spaces $\gX$ and $\gY$ and a matrix representation with respect to a particular basis are not distinguished notationally. The image and the kernel of matrix $A$ are denoted by $\ima\,A$ and $\ker\,A$, respectively. The Moore-Penrose pseudo-inverse of $A$ is denoted by $A^\dagger$. When $A$ is square, we denote by $\sigma(A)$ the spectrum of $A$. If $\gJ \subseteq \gX$, the restriction of the map $A$ to $\gJ$ is denoted by $A\,|\gJ$. If $\gX=\gY$ and $\gJ$ is $A$-invariant, the eigenstructure of $A$ restricted to $\gJ$ is denoted by $\sigma\,(A\,| \gJ)$. If $\gJ_1$ and $\gJ_2$ are $A$-invariant subspaces and $\gJ_1\,{\subseteq}\,\gJ_2$, the mapping induced by $A$ on the quotient space $\gJ_2 / \gJ_1$ is denoted by $A\,| {\gJ_2}/{\gJ_1}$, and its spectrum is denoted by $\sigma\,(A\,| {\gJ_2}/{\gJ_1})$.
The symbol $\oplus$ stands for the direct sum of subspaces.  The symbol $\uplus$ denotes union with any common elements repeated.
Given a map $A: \gX \longrightarrow \gX$ and a subspace $\gB$ of $\gX$, we denote by $\langle A, \gB \rangle$ the smallest $A$-invariant subspace of $\gX$ containing $\gB$. 
Given a complex matrix $M$, the symbol ${M}^\ast$ denotes the conjugate transpose of $M$. Moreover, we denote by $M_i$ its $i$-th row and by $M^j$ its $j$-th column, respectively.
 Given a finite set $S$, the symbol $2^S$ denotes the power set of $S$, while ${\rm {card}}(S)$ stands for the cardinality of $S$.
We recall that a subset of $\real^n$ is a {\em Zariski open set}  if it is nonempty and its complement is formed by the solutions to finitely many polynomial equations $\phi_i(p_1,\ldots,p_n)=0$ ($i \in \{1,\ldots,k\}$) where the coefficients of the polynomials $\phi_i$ are real. Let $\Pi$ be a property defined on $\real^n$ (i.e., $\Pi: \real^{n} \longrightarrow \{0,1\}$). Consider the set $\Phi=\{\pi\in \real^n\,|\,\Pi(\pi)=1\}$. If $\Phi$ is a Zariski open set, then we say that {\em almost all} $\pi\in \real^n$ satisfy $\Pi$, or that $\Pi$ is {\em generic} in $\real^n$.

\section{Problem Formulation}
\label{PF}
In what follows, whether the underlying system evolves in continuous or discrete time makes only minor differences and, accordingly, the time index set of any signal is denoted by $\mathbb{T}$, on the understanding that this represents either $\real^+$ in the continuous time or $\mathbb{N}$ in the discrete time. The symbol $\complex_g$ denotes either the open left-half complex plane $\complex^-$ in the continuous time or the open unit disc $\complex^\circ$ in the discrete time. A matrix $M \in \real^{n \times n}$ is said to be {\em asymptotically stable} if $\sigma(M) \subset \complex_g$. Finally, we say that $\l \in \complex$ is stable if $\l \in \complex_g$.
Consider the LTI system $\Sigma$ governed by
 \bea
 \Sigma: \
 \left\{ \begin{array}{lcr}
 \gD\,x(t) \ns&\ns = \ns&\ns A\,x(t)+B\,u(t),\;\;\;\; x(0)=x_0, \label{syseq1}\hfill\cr
 y(t) \ns&\ns = \ns&\ns C\,x(t)+D\,u(t),\hfill \end{array} \right.
 \label{sys}
 \eea
where, for all $t \in \mathbb{T}$, $x(t) \in \gX=\real^n$ is the
state, $u(t) \in \gU=\real^m$ is the control input, $y(t) \in \gY=\real^p$ is the output, and $A$, $B$, $C$ and $D$ are appropriate
dimensional constant matrices. The operator
$\gD$ denotes either the time derivative in the continuous time,
i.e., $\gD\, x(t)=\dot{x}(t)$, or the unit time shift
 in the discrete time, i.e., $\gD\, x(t)=x(t+1)$. 
 We assume with no loss of generality that all the columns of $\left[ \begin{smallmatrix} B \\[1mm] D \end{smallmatrix}\right]$ are linearly independent. We also assume that all the rows of
 $[\, C \;\; D \,]$ are linearly independent. 
We recall that the Rosenbrock system matrix is defined as the matrix pencil
\bea
\label{ros}
P_{\scriptscriptstyle \Sigma}(\lambda) \defi \bmat{cc} A-\lambda\,I_n & B \\ C & D \emat
\eea
in the indeterminate $\lambda \in \complex$, see e.g. \cite{Rosenbrock-70}. The invariant zeros of $\Sigma$ are the values of $\lambda \in \complex$ for which the rank of $P_{\scriptscriptstyle \Sigma}(\lambda)$ is strictly smaller than its normal rank.\footnote{The normal rank of a rational matrix $M(\lambda)$ is defined as
$\normrank M(\lambda) \defi \displaystyle\max_{\lambda \in \complex} \rank M(\lambda)$. The rank of $M(\lambda)$ is equal to its normal rank for all but finitely many $\lambda \in \complex$.} More precisely,
 the invariant zeros are the roots of the non-zero polynomials on the principal diagonal of the Smith form of $P_{\scriptscriptstyle \Sigma}(\lambda)$, see e.g. \cite{Trentelman-SH-01}.
Given an invariant zero $\lambda=z \in \complex$, the rank deficiency of $P_{\scriptscriptstyle \Sigma}(\lambda)$ at the value $\lambda=z$ is the geometric multiplicity of the invariant zero $z$, and is equal to the number of elementary divisors (invariant polynomials) of $P_{\scriptscriptstyle \Sigma}(\lambda)$ associated with the complex frequency $\lambda=z$. The degree of the product of the elementary divisors of $P_{\scriptscriptstyle \Sigma}(\lambda)$ corresponding to the invariant zero $z$ is the algebraic multiplicity of $z$, see \cite{Trentelman-SH-01}. The set of invariant zeros of $\Sigma$ is denoted by $\gZ$, and the set of minimum-phase invariant zeros of $\Sigma$ is $\gZ_g\defi \gZ\cap \complex_g$.

We denote by $\gV^\star$ the largest output-nulling subspace of $\Sigma$, i.e., the largest subspace $\gV$ of $\gX$ for which a matrix $F\,{\in}\,\mathbb{R}^{m\,{\times}\,n}$ exists such that $(A+B\,F)\,\gV\subseteq \gV \subseteq \ker (C+D\,F)$. Any real matrix $F$ satisfying this inclusion is called a {\it friend \/} of $\gV$. We denote by $\mathfrak{F}(\gV)$ the set of friends of $\gV$.
 The symbol $\gR^\star$ denotes the so-called {\em output-nulling reachability subspace} on $\gV^\star$, and is the smallest $(A\,{+}\,B\,F)$-invariant subspace of $\gX$ containing $\gV^\star\,{\cap}\,B\,\ker\,D$, where $F\,{\in}\,\mathfrak{F}(\gV^\star)$.
 The closed-loop spectrum can be partitioned as
$\sigma(A+B\,F)=\sigma(A+B\,F\,|\,\gV^\star)\uplus \sigma(A+B\,F\,|\,\gX/\gV^\star)$, where
 $\sigma(A+B\,F\,|\,\gV^\star)$ is the spectrum of $A+B\,F$ restricted to $\gV^\star$, and
$\sigma(A+B\,F\,|\,\gX/\gV^\star)$ is the spectrum of the mapping induced by $A+B\,F$ on the quotient space $\gX/\gV^\star$. 
 The eigenstructure of $A+B\,F$ restricted to $\gV^\star$ can be further split into two disjoint parts: the eigenstructure $\sigma(A+B\,F |\gR^\star)$ is freely assignable\footnote{An assignable set of eigenvalues is always intended to be a set of complex numbers mirrored with respect to the real axis. }
 with a suitable choice of $F$ in $\mathfrak{F}(\gV^\star)$. The eigenstructure  $\sigma\,(A+B\,F | {\gV^\star}/{\gR^\star})$ is fixed for all the choices of $F$ in $\mathfrak{F}(\gV^\star)$ {and coincide with the invariant zero structure of $\Sigma$, \cite[Theorem 7.19]{Trentelman-SH-01}}.
 Finally, we use the symbol $\gV^\star_g$ to denote the largest stabilisability output-nulling subspace of $\Sigma$, i.e., the largest subspace $\gV$ of $\gX$ for which a matrix $F\,{\in}\,\mathbb{R}^{m\,{\times}\,n}$ exists such that $(A+B\,F)\,\gV\subseteq \gV \subseteq \ker (C+D\,F)$ and $\sigma(A+B\,F\,|\,\gV) \subset \complex_g$. There holds $\gR^\star \subseteq \gV^\star_g \subseteq \gV^\star$. \\[-2mm]

 \subsection{The tracking control problem}

In this paper, we are concerned with the problem of the design of a
state-feedback control law for (\ref{sys}) such that for all initial conditions the output $y$ tracks a step reference $r \in \gY$ with zero steady-state error and is monotonic in all components.
If $y$ asymptotically tracks the constant reference $r$ and is monotonic, then it is also both non-overshooting and non-undershooting. The converse is obviously not true in general.
 The following standing assumption is standard for tracking control problems (see e.g. \cite{He-CW-05}), and ensures that any given constant
reference target $r$ can be tracked from any initial
condition $x_0 \in \gX$:\\[-3mm]
 \begin{assumption}
 \label{Ass1}
System $\Sigma$ is right invertible and stabilisable, and
 $\Sigma$ has no invariant
zeros at the origin in the continuous time case, or at $1$ in the discrete case.\\[-3mm]
\end{assumption}

Assumption \ref{Ass1} generically hold when $m \ge p$.
Under Assumption \ref{Ass1}, the standard method for designing a tracking controller for a
step reference signal is carried out as follows. Given the step reference $r \in \gY$ to track, choose a feedback gain matrix
$F$ such that $A + B\,F$ is asymptotically stable: this is always possible since the pair $(A,B)$ is assumed to be stabilisable. Let us then choose two vectors $x_{\rm ss} \in
\gX$ and $u_{\rm ss} \in \gU$ that, for the given $r \in \gY$, satisfy 
\begin{eqnarray}
 \label{k2}
\left\{ \begin{array}{rcl}
0 \ns&\ns = \ns&\ns A\,x_{\rm ss}+B\,u_{\rm ss} \\
r \ns&\ns = \ns&\ns C\,x_{\rm ss}+D\,u_{\rm ss} \end{array} \right. \quad \text{or} \quad 
\left\{ \begin{array}{rcl}
 x_{\rm ss} \ns&\ns = \ns&\ns A\,x_{\rm ss}+B\,u_{\rm ss} \\ 
  r \ns&\ns = \ns&\ns C\,x_{\rm ss}+D\,u_{\rm ss}
\end{array} \right.
\end{eqnarray}
in the continuous and in the discrete case, respectively. 
 Such pair of vectors $x_{\rm ss} \in \gX$ and $u_{\rm ss} \in \gU$ exist since {\em (i)} right invertibility ensures that the system matrix pencil
$P_{\scriptscriptstyle \Sigma}(\lambda)$ is of full row-rank for all but finitely many $\lambda \in \complex$, see \cite[Theorem 8.13]{Trentelman-SH-01}, and, as already recalled, the values $\lambda \in \complex$ for which $P_{\scriptscriptstyle \Sigma}(\lambda)$ loses rank are invariant zeros of $\Sigma$; {\em (ii)} in the continuous (resp. discrete) time case, the absence of invariant zeros at the origin (resp. at $1$) guarantees that the matrix $P_{\scriptscriptstyle \Sigma}(0)$ (resp. $P_{\scriptscriptstyle \Sigma}(1)$) is of full row-rank. As such, Assumption \ref{Ass1} guarantees that the linear system (\ref{k2}) is always solvable in $\left[ \begin{smallmatrix} x_{\rm ss} \\[1mm] u_{\rm ss} \end{smallmatrix} \right]$. Now, applying the state-feedback control law
 \be
 u(t) = F\,\Big( x(t)-x_{\rm ss} \Big) + u_{\rm ss}
 \label{ulaw}
 \ee
 to (\ref{sys}) and using the change of variable $\xi \defi x -x_{\rm ss}$ gives the closed-loop
 autonomous system
 \be
 \Sigma_{aut}: \
\left\{ \begin{array}{lcr}
\gD\,{\xi}(t) \ns&\ns = \ns&\ns (A+B\,F)\,\xi(t),\;\;\quad \xi(0)=x_0-x_{\rm ss},\hfill\cr
 y(t) \ns&\ns = \ns&\ns (C+D\,F)\,\xi(t)+ r. \hfill \end{array} \right.
 \label{syschom}
 \ee
Since $A+B\,F$ is asymptotically stable, $x$ converges to $x_{\rm
ss}$, $\xi$ converges to zero and $y$ converges to $r$ as $t$ goes to infinity. We shall refer to $\xi$ as the {\em error state coordinates}.

 \subsection{Achieving a globally monotonic response with any desired convergence rate}

In this paper we are concerned with the general problem of finding a gain matrix $F$ such that the closed-loop system obtained using \eqref{ulaw} in \eqref{sys} achieves a monotonic response at any desired rate of convergence, from all initial conditions. We shall describe this property as {\em global monotonicity}. We describe the problem as follows in terms of the tracking error \[
\epsilon(t)\defi y(t)-r(t) = \bsmat \epsilon_1(t) \\[-1mm] \vdots \\[1mm] \epsilon_p(t)\esmat.
\]

\begin{problem}
\label{prob:mono0}
Let $\rho \in \real$, such that $\rho<0$ in the continuous time and $\rho \in [0,1)$ in the discrete time. Find a state-feedback matrix $F$ such that applying (\ref{ulaw}) with this $F$ to $\Sigma$ yields an asymptotically stable closed-loop system $\Sigma_{aut}$ for which the tracking error term $\epsilon(t)$ converges monotonically to 0 at a rate at least $\rho$ in all outputs, from all initial conditions. Specifically, we require that in the continuous time
\begin{equation}
 \forall \xi(0) \in \gX,\	\forall k \in \{1, \dots, p\}, \ \exists \,\beta_k \in \mathbb{R} \, : \;\; | \epsilon_k(t) | \leq \beta_k \, \exp (\rho\, t) 
\label{eq:error2}
\end{equation}
for all $t \in \real_+$, where $\epsilon_k(t)$ is strictly monotonic in $t$, and in the discrete time
\begin{equation}
 \forall \xi(0) \in \gX,\	\forall k \in \{1, \dots, p\}, \ \exists \,\beta_k \in \mathbb{R} \, : \;\; | \epsilon_k(t) | \leq \beta_k \, {\rho^{t}} %
\label{eq:error2d}
\end{equation}
for all $t \in \mathbb{N}$, where, again, $\epsilon_k(t)$ is strictly monotonic in $t$.
\end{problem}

If we are able to obtain a tracking error $\epsilon(t)$ that consists of a single exponential per component in the continuous time or a single power per component in the discrete time, i.e., 
\begin{equation} 
\label{eq:error}
\epsilon(t) =
 \begin{bmatrix} \beta_1\,\exp(\lambda_1\,t)\\[-2mm]\vdots \\[-1mm] \beta_p\,\exp(\lambda_p\,t) \end{bmatrix} \quad \text{or} \quad
\epsilon(t) =
 \begin{bmatrix} \beta_1\,\lambda_1^{t} \\[-2mm]\vdots \\[-1mm] \beta_p\,\,\lambda_p^{t} \end{bmatrix},
 \end{equation}
 respectively, and we can choose each $\lambda_k$ in such a way that $\lambda_k \leq \rho<0$ in the continuous time and $0 \le \lambda_k \leq \rho$ in the discrete time, 
then we solve Problem \ref{prob:mono0}. Indeed, asymptotically stable exponentials of $\lambda_k$ or powers of $\lambda_k$ are monotonic functions. This suggests that a possible way of solving Problem \ref{prob:mono0} consists in the solution of the following problem.

\begin{problem}
\label{prob:mono1}
Let $\rho \in \real$, such that $\rho<0$ in the continuous time and $\rho \in [0,1)$ in the discrete time. Find a feedback matrix $F$ such that applying (\ref{ulaw}) with this $F$ to $\Sigma$ yields an asymptotically stable closed-loop system $\Sigma_{aut}$ for which, from all initial conditions, the tracking error term is given by
(\ref{eq:error})
for some real coefficients $\{\beta_k \}_{k=1}^p$ depending only upon $\xi(0)$ and for some real values $\lambda_1, \lambda_2, \cdots, \lambda_p$ 
satisfying in the continuous time $\lambda_k \leq \rho<0$ and in the discrete time $0 \le \lambda_k \leq \rho$. 
\end{problem}

Clearly, solutions of Problem \ref{prob:mono1} also solve Problem \ref{prob:mono0}.
However, the following result shows that the converse is also true: the only way to obtain a feedback matrix ensuring global monotonic tracking is to obtain a tracking error as in (\ref{eq:error}).\\[-3mm]

\begin{lemma} 
\label{lem:mono}
Problem~\ref{prob:mono0} is equivalent to Problem~\ref{prob:mono1}.
\end{lemma}
\proof
Let us consider the continuous time case, the discrete case being entirely equivalent.
 Let $\rho<0$. If $F$ and $\lambda_1, \lambda_2, \cdots, \lambda_p$ solve Problem~\ref{prob:mono1} with respect to $\rho$, then the outputs $\epsilon_k(t)$ satisfy (\ref{eq:error}), and hence also (\ref{eq:error2}) or (\ref{eq:error2d}). 
 Next, assume that the feedback matrix $F$ solves Problem~\ref{prob:mono0} for a certain $\rho \in \real^-$. The tracking error is the output of the autonomous system $\Sigma_{aut}$, and its components can be written as 
$\epsilon_k(t)=[C+D\,F]_k\,e^{(A+B\,F)\,t}\,\xi_0$,
 where $[C+D\,F]_k$ denotes the $k$-th row of $C+D\,F$. Let us change coordinates, and let us write the pair $(C+D\,F,A+B\,F)$ in the standard observability form, in which the observable part (whose dimension is denoted by $l$) is in turn expressed in the observability canonical form, i.e.,
 \beann
 C+D\,F \ns&\ns = \ns&\ns [\begin{array}{cccc|ccc} \! 0  \! \!  &  \! \!  \ldots  \!  \! & \!   \! 0 \!   \! &  \! \!  1  \!  \! &  \! \!  0  \!  \! &  \! \!  \ldots \!  \!  &  \! \!  0 \!  \end{array}] \\ 
 A+B\,F \ns&\ns = \ns&\ns \bmat{c|c} \begin{smallmatrix} 
\!\!\!\! 0 && &&  && -\alpha_0 \!\!\!\! \\[-2mm]
\!\!\!\!  1 && {\tiny \ddots} && && -\alpha_1 \!\!\!\! \\[-2mm]
\!\!\!\!  && {\tiny \ddots} && 0 && \!\!\!\! \\[1mm]
\!\!\!\!  &&&& 1 && -\alpha_{l-1} \!\!\!\! \end{smallmatrix} & 0 \\[1mm] \hline Z_{21} & Z_{22} \emat
 \eeann
 and, accordingly, $\xi_0$ in the new coordinates is written as
 $\xi_0= [\begin{array}{ccc|ccc} \xi_{0,1} & \ldots & \xi_{0,l} & \xi_{0,l+1} & \ldots \xi_{0,n} \end{array}]^\top$. The Laplace transform of $\epsilon_k(t)$ is therefore
 \[
 \gL[\epsilon_k]=\left(   \! [\begin{array}{cccc}  \! \! 1 \! \!   & \!  s   \! \! & \!  \ldots  \! \!  & \!  s^{l-1}  \!  \! \!  \end{array} ]  \!  \! \bsmat \xi_{0,1} \\[1mm]
 \xi_{0,2} \\[-1.5mm]
\vdots \\
\xi_{0,l} \esmat\right) \frac{1}{s^l+\alpha_{l-1}\,s^{l-1}+ \ldots + \alpha_0},
\]
which shows that the numerator of $\gL[\epsilon_k]$ { is a polynomial, say $N(s)$, and the denominator, say $D(s)$, is the characteristic polynomial of the observable part of $(C+D\,F,A+B\,F)$.
By suitably choosing the initial condition $\xi_0$,  the degree $d$ of $N(s)$ can be selected arbitrarily in the range  $\{0,1,2\dots,l-1\}$ and the coefficients of $N(s)$ can be selected arbitrarily in $\real^{d+1}$.
Thus $\epsilon_k(t)$ is a linear combination of the modes of the observable subsystem} with arbitrary coefficients, and can be written as
 \bea
\epsilon_k(t) \ns&\ns = \ns&\ns \sum_{i=1}^{\rho} \sum_{j=1}^{m_i}  \tilde{\beta}_{k,i}\, t^{j-1} e^{\lambda_i t}+\sum_{i=1}^{c} \sum_{j=1}^{\tilde{m}_i} [\hat{\beta}_{k,i}^\prime\, t^{j-1} \,e^{\sigma_i t}\,\cos(\omega_i\,t) \nonumber \\
\ns&\ns \ns&\ns +\hat{\beta}_{k,i}^{\prime \prime}\, t^{j-1} \,e^{\sigma_i t}\,\sin(\omega_i\,t)], \label{error}
\eea
where $\lambda_1,\ldots,\lambda_{\rho}$ are the real eigenvalues of the observable subsystem with associated algebraic multiplicities $m_1,\ldots, m_{\nu}$ and where $\mu_1,\ldots,\mu_c,\overline{\mu}_1,\ldots,\overline{\mu}_c$ are the complex eigenvalues of the observable subsystem, and the algebraic multiplicities associated with $\mu_1,\ldots,\mu_c$ are $\tilde{m}_1,\ldots, \tilde{m}_{c}$, where $\sigma_i = \mathfrak{Re}\{{\mu}_i\}$ and $\omega_i = \mathfrak{Im}\{{\mu}_i\}$. Finally, 
 the real coefficients $\tilde{\beta}_{k,i}$, $\hat{\beta}_{k,i}^\prime$ and $\hat{\beta}_{k,i}^{\prime \prime}$  can be made arbitrary by choosing suitable initial conditions.
  In particular, we can pick an arbitrary one of the $l$ modes appearing in (\ref{error})
and select the initial conditions in such a way that $\epsilon_k(t)$ is that mode.
Since by assumption the response is monotonic from any initial condition, only components of the form $\exp({\lambda}_i \,t)$ for real ${\lambda}_i$ can appear in each $\epsilon_k$, because for any real ${\lambda}_i<0$, the function $t^{j-1}\,\exp({\lambda}_i\, t)$ is monotonic only if $j =1$ while the functions $t^{j-1} \exp(\sigma_i\, t)\cos(\omega_i\, t)$ and $t^{j-1} \exp(\sigma_i \,t)\sin(\omega_i\, t)$ are monotonic only if $j =1$ and $\omega_i=0$. 
In other words, the polynomial $D(s)$ can have only real, simple, negative roots.
Thus for each output $\epsilon_k(t)$ we have 
$\epsilon_k(t) = \sum_{i=1}^{\nu} \tilde{\beta}_{k,i} \exp({\lambda}_i t)$,
where the real coefficients $\tilde{\beta}_{k,i}$ can be made arbitrary by choosing suitable initial conditions. We now use again the fact that each response is monotonic from all initial conditions. From Lemma A.1 of \cite{Schmid-N-09}, if $\epsilon_k(t)$ is a linear combination of two or more negative real exponential functions, it will change sign (and hence not be monotonic) for some values of the coefficients $\tilde{\beta}_{k,i}$. Thus, for each $k \in \{1,\ldots,p\}$ we must have
$\epsilon_k(t) = \tilde{\beta}_{k,i} \exp({\lambda}_i t)$
for some eigenvalue ${\lambda}_i$ and some real coefficient $\beta_k \defi \tilde{\beta}_{k,i}$. Thus, Problem 2 is solved. 
\endproof

\ \\[-6mm]

For conciseness, let us define the set $\Lambda_g$ to be equal to $\real^-$ in the continuous time and $[0,1)$ in the discrete time.

Another important and useful problem is one in which the requirements include a specified choice of the closed-loop modes that are visible in each component of the tracking error:
\ \\[-4mm]
\begin{problem} 
\label{prob:mono2}
Let $\lambda_1, \cdots, \lambda_p \in \Lambda_g$. Find a feedback matrix $F$ such that applying (\ref{ulaw}) 
to $\Sigma$ yields an asymptotically stable closed-loop system $\Sigma_{aut}$ for which, from all initial conditions and all step references, the tracking error term is given by (\ref{eq:error}). 
\end{problem}

\section{Global monotonicity: the intuitive idea} 
\label{sec:intu}
In the previous section, we observed that in order for the problem of global monotonic tracking to be solvable, we need to distribute (at most) $p$ modes evenly into the tracking error with one mode per error component. To achieve this goal, all the remaining closed-loop modes have to be made invisible from the tracking error. 
 If this is possible, then the step response is guaranteed to be monotonic for any initial condition, and therefore also non-overshooting and non-undershooting. 
The converse is also true, as shown in Lemma \ref{lem:mono}.
If we are able to render more than $n-p$ modes invisible at $\epsilon(t)$, one or more components of the tracking error can be rendered identically zero, and  for those components instantaneous tracking may also be achieved, in which the output component immediately takes the desired reference value.
Our aim is to find conditions under which a gain matrix $F$ can be obtained to deliver the single mode structure (\ref{eq:error}) for any initial condition. 
Consider $\Sigma_{aut}$ in (\ref{syschom}), which can be re-written as
 \be
 \Sigma_{aut}: \
\left\{ \begin{array}{lcr}
\gD\,{\xi}(t) \ns&\ns = \ns&\ns A\,\xi(t)+B\,\omega(t),\;\;\; \hfill\cr
\epsilon(t) \ns&\ns = \ns&\ns C\,\xi(t)+ D\,\omega(t), \hfill \end{array} \right.
 \label{syschom1}
 \ee
where $\omega(t)=F\,\xi(t)$. Clearly, $\Sigma_{aut}$ can be identified with the quadruple $(A,B,C,D)$. 
The task is now to 
find a feedback matrix $F$ such that the new control $\omega(t)=F\,\xi(t)$ guarantees that for every initial condition $\xi_0\in\gX$ the tracking error $\epsilon(t)$ is characterised by a single stable real mode per component.
Let $j \in \{1,\ldots,p\}$. Let $\lambda_j \in \Lambda_g$. Consider a solution $v_j$ and $w_j$ of the linear equation
\bea
\label{ML0}
\left[ \begin{array}{ccc} \! A-\lambda_j\,I_n \! & \! B \! \\ \! C \! & \! D \! \end{array} \right]\left[ \begin{array}{ccc} v_j \\ w_j \end{array} \right]=\left[ \begin{array}{ccc} 0 \\ \beta_j\,e_j \end{array} \right],
\eea
where $\beta_j \neq 0$ and $e_j$ is the $j$-th vector of the canonical basis of $\gY$. Notice that (\ref{ML0}) always has a solution in view of the right-invertibility of $\Sigma$. By choosing $F$ such that $F\,v_j=w_j$, we find $(A+B\,F)\,v_j=\lambda_j\,v_j$ and $(C+D\,F)\,v_j=\beta_j\,e_j$. Hence, from (\ref{error}) we know that for any initial error state $\xi_0\in\spanR\{v_j\}$ the response associated with the control $\omega(t)=F\,\xi(t)$ is
\bea
\label{newalpha}
\epsilon(t)=\left[ \begin{array}{ccc} 0 \\[-2mm] \tiny{\vdots} \\[-0mm] \gamma_j\,\exp(\lambda_j\,t) \\[-1mm] \tiny{\vdots} \\[-0.5mm] 0 \end{array} \right] \left. \begin{array}{ccc} \phantom{0} \\[-0mm] \phantom{\vdots} \\[-0mm] \leftarrow j \phantom{q\,\exp(\lambda\,t)} \\[-1mm] \phantom{\vdots} \\[-0mm] \phantom{0} \end{array} \right.
\eea
where $\gamma_j$ depends on the particular initial state $\xi_0$. Considering $\l_1,\ldots,\l_p \in \Lambda_g$, 
by applying this argument for all components of the tracking error, we obtain a set of solutions $\left[\begin{smallmatrix} v_1 \\[1mm] w_1 \end{smallmatrix} \right], \left[\begin{smallmatrix} v_2 \\[1mm] w_2 \end{smallmatrix} \right], \ldots,\left[\begin{smallmatrix} v_p \\[1mm] w_p \end{smallmatrix} \right]$ of (\ref{ML0}). If $v_1,\ldots,v_p$ are linearly independent, we can choose
 $F$ to be such that $F\,v_i=w_i$ for all $i \in \{1,\ldots,p\}$. Then, for every $\xi_0 \in \spanR \{v_1,v_2,\ldots,v_p\}$, by superposition we find
\bea
\epsilon(t) \ns&\ns \!\!=\!\! \ns&\ns  \! \left[ \begin{array}{ccc} \!\!\!\!\!  \gamma_1\,e^{\lambda_1\,t} \!\!\!\!\!   \\[-0mm]  \!\!\!\!\!  0  \!\!\!\!\!  \\[-2mm] \!\!\!\!\!    \tiny{\vdots}  \!\!\!\!\!   \\[-0mm]  \!\!\!\!\!  0 \!\!\!\!\!   \end{array} \right] \!\!+\!\! \left[ \begin{array}{ccc} \!\!\!\!\! 0 \!\!\!\!\!   \\[-0mm]  \!\!\!\!\!   \gamma_2\,e^{\lambda_2\,t}  \!\!\!\!\!  \\[-2mm] \!\!\!\!\!   \tiny{\vdots}  \!\!\!\!\!   \\[-0mm]  \!\!\!\!\!  0 \!\!\!\!\!   \end{array} \right]\!\!+\ldots+\!\!
\left[ \begin{array}{ccc} 0 \\[-0mm] 0 \\[-2mm]\ \tiny{\vdots}  \\[-0mm]  \!\!\!\!\!  \gamma_p\,e^{\lambda_p\,t} \!\!\!\!\!   \end{array} \right] \!=\!
\left[ \begin{array}{ccc} \!\!\!\!\!   \gamma_1\,e^{\lambda_1\,t}   \!\!\!\!\!  \\[-0mm]  \!\!\!\!\!  \gamma_2\,e^{\lambda_2\,t}  \!\!\!\!\!   \\[-2mm] \tiny{\vdots} \\[-0mm] \!\!\!\!\!   \gamma_p\,e^{\lambda_p\,t}  \!\!\!\!\!   \end{array} \right]\!\!. \phantom{pipi}
\label{risp}
\eea
However, this result only holds when $\xi_0 \in \spanR \{v_1,v_2,\ldots,v_p\}$.
In order for this response to be achievable from any initial condition, we also need to render  the remaining $n-p$ closed-loop modes invisible at $\epsilon(t)$. This task can be accomplished by exploiting the supremal stabilisability output-nulling subspace $\gV^\star_g$ of the system, which is defined as the largest subspace of $\gX$ for which a friend $F$ exists such that, for every initial state lying on it, the corresponding state feedback generates a state trajectory that asymptotically converges to zero while the corresponding output (the tracking error in the present case) remains  at zero.
We shall see in Section \ref{background} that a basis for $\gV^\star_g$ can always be obtained 
as the image of a matrix $[\begin{array}{cccc} V_1 & V_2 & \ldots & V_d\end{array}]$ that satisfies
\bea
\label{RosVg}
\left[ \begin{array}{ccc} \! A-\mu_j\,I_n \! & \! B \! \\ \! C \! & \! D \! \end{array} \right]\left[ \begin{array}{ccc} V_{j} \\ W_{j} \end{array} \right]=0,
\eea
for some other matrix $[\begin{array}{cccc} W_1 & W_2 & \ldots & W_d\end{array}]$ partitioned conformably, 
where $\{\mu_1,\ldots,\mu_t\}$ are the (distinct) minimum-phase invariant zeros of $\Sigma$ and $\{\mu_{t+1},\ldots,\mu_d\}$ are arbitrary and stable (let us assume for the moment that they are real and distinct for simplicity).
If the dimension of $\gV^\star_g+\spanR \{v_1,v_2,\ldots,v_p\}$ is equal to $n$, every initial state $\xi_0\in\gX$ can be decomposed as the sum $\xi_v+\xi_r$, where $\xi_v\in \gV^\star_g$ and $\xi_r \in \spanR \{v_1,v_2,\ldots,v_p\}$. 
 If for the sake of argument we have $\dim \gV^\star_g=n-p$, 
 and we can find a set of linearly independent columns $\{v_{p+1},\ldots,v_n\}$ from the columns of $[\begin{array}{cccc} \!\!\! V_1\!\!\! & V_2 \!\!\!& \ldots\!\!\! & V_d\!\!\!\end{array}]$ that is linearly independent from $\{v_1,\ldots,v_p\}$, we can take $w_{p+1},\ldots,w_{n}$ to be the columns of $W_g$ that correspond to $v_{p+1},\ldots,v_{n}$, and construct the feedback control $\omega(t)=F\,\xi(t)$ where $F$ is such that
 $F\,[\begin{array}{cccccc} \!\! v_1\!\!\!\! &\!\!\ldots \!\!& \!\!v_p \!\!& \!\!v_{p+1}\!\! &\!\!\ldots \!\!&\!\!v_{n}\!\!\end{array}]=[\begin{array}{cccccc} \!\!w_1\!\! &\!\!\ldots\!\! &\!\! w_p\!\! & \!\!w_{p+1}\!\!& \!\!\ldots\!\!& \!\!w_{n}\!\!\end{array}]$, the response associated with $\xi_v$ is identically zero, while the one associated with $\xi_r$ is still given by (\ref{risp}). Hence, the tracking error can be written as in (\ref{eq:error}) for any $\xi_0 \in \gX$. The closed-loop eigenvalues obtained with  $F$ are given by the union of $\{\l_1,\ldots,\l_p\}$, with the set of $\mu_j$ that are associated with the columns $\{v_{p+1},\ldots,v_n\}$ chosen from $[\begin{array}{cccc} \!\!\! V_1\!\!\! & V_2 \!\!\!& \ldots\!\!\! & V_d\!\!\!\end{array}]$.

We now introduce a running example that embodies all those system theoretic characteristics that are perceived as the major difficulties in achieving monotonic tracking. This system is MIMO, non-strictly proper, uncontrollable (but obviously stabilisable) and is characterised by 3 non-minimum phase zeros. To the best of the authors' knowledge, there are no methods available in the literature that can solve the tracking problem for MIMO systems with a guaranteed monotonic response under such assumptions, and in particular in the presence of three non-minimum phase invariant zeros. 
We also want to stress that this problem is solved here in closed form.

\begin{example}
\label{exe0}
{Consider the non-strictly proper continuous-time LTI system $\Sigma$ in (\ref{sys}) with
\beann
A \ns&\ns = \ns&\ns \bsmat
 -6 && 0 && 0 && 0 && 0\\[1mm]
 3 && 3 && 0 && 0 && 0\\[1mm]
 0 && 0 && 2 && 0 && 2\\[1mm]
 -1 && 0 && 2 && 0 && 0\\[1mm]
 -2 && 0 && 0 && 0 && 2\esmat\!, \quad B=\bsmat
 0 && 0 && 0 && 0\\[1mm]
 0 && 0 && 0 && -3\\[1mm]
 0 && 4 && 2 && 0\\[1mm]
 1 && -1 && 0 && -1\\[1mm]
 0 && -1 && 0 && 0 \esmat\!,\\
 C \ns&\ns = \ns&\ns \bsmat
 -1 && 0 && 0 && 0 && 0\\[1mm]
 3 && 0 && 0 && 0 && 9\\[1mm]
 1 && 0 && 0 && 0 && 0\esmat\!, \qquad D=\bsmat
 0 && 0 && -2 && 0\\[1mm]
 0 && 3 && -3 && -3\\[1mm]
 0 && 0 && 2 && -2 \esmat.
 \eeann
 We want to find a feedback matrix $F$ such that the output of $\Sigma$ monotonically tracks a unit step in all output components, and the assignable closed-loop eigenvalues are equal to $\lambda_1=\lambda_2=\lambda_3=-1$ for the corresponding error components.
 This system is right invertible but not left invertible. It is not reachable but it is stabilisable, since the only uncontrollable eigenvalue is $-6$, and is equal to the only minimum-phase invariant zero $z_1=-6$.  The other three invariant zeros of the system are non-minimum-phase, and their values are $z_2=2$, $z_3=3$ and $z_4=5$.
A feedback matrix which solves the problem is 
 \bea
 \label{f1}
 F=\bsmat
  \frac{925}{198}  &  \frac{4}{3}  &  -2  &  -1  &  3   \\[1mm]
  -\frac{39}{22}  &  0  &  0  &  0  &  3   \\[1mm]
  \frac{107}{88}  &  0  &  -\frac{3}{2}  &  0 &  -7  \\[1mm]
  \frac{4}{9} &  \frac{4}{3}  &  0  &  0  &  0  \esmat.
 \eea
  The closed-loop eigenvalues are $\sigma(A+B\,F)=\{-1,-6\}$ where the multiplicity of the eigenvalue $-1$ is four. If the reference is $r=[\begin{array}{cccc} \! 2 \! & \! 2 \! & \! 2 \! \end{array}]^\top$, using (\ref{k2}) we obtain $x_{\rm ss} =[\begin{array}{ccccc} \!0 \! & \! -2  \! & \! {10}/{3} \! & \! 0  \! & \!  -{7}/{15}\!\end{array}]^\top$ and $u_{\rm ss}=[\begin{array}{cccc} \!-{48}/{5}  \! & \! 
 -{14}/{15} \! & \!  -1 \! & \! -2 \!\end{array}]^\top$.
 Given an arbitrary initial condition $\xi_0 \in \gX$, the tracking error that follows from the application of the control law $u(t) = F\,( x(t)-x_{\rm ss} ) + u_{\rm ss}$ with the feedback matrix $F$ in (\ref{f1}), yields 
 $\epsilon(t)=[\, \gamma_1\,e^{-t} \;\; \gamma_2\,e^{-t} \;\; \gamma_3\,e^{-t} \,]^\top$, which has the single mode form of (\ref{risp}). Thus, the system exhibits a globally monotonic step response.
 \endex
 }
 \end{example}

Since the solvability condition for global monotonicity is given in terms of the dimension of the subspace $\gV^\star_g+\spanR \{v_1,\ldots,v_p\}$, and $v_1,\ldots,v_p$ depend on the choice of the closed-loop eigenvalues $\lambda_1,\ldots,\lambda_p$, the solvability condition seems to depend on the particular choice of the closed-loop eigenvalues. 
The question at this point is: how does the choice of the closed-loop eigenvalues affect the dimension of $\gV^\star_g+\spanR \{v_1,\ldots,v_p\}$? Are there good and bad choices of the closed-loop eigenvalues? More generally, can we find alternative solvability conditions given solely in terms of the system structure and not in terms of a choice of eigenvalues? These are the crucial points that will be addressed in the sequel.

\section{Mathematical background}
\label{background}
In the previous section, it has been shown that it is necessary to render at least $n-p$ modes invisible at the tracking error in order to achieve global monotonicity. We now provide the basic tools which will be used for such a purpose.
The first of these tools is the subspace $\gV^\star_g$,
which is made up of the sum of two parts. The first is the subspace $\gR^\star$. The second is, loosely, the subspace spanned by the directions associated with the minimum-phase invariant zeros of $\Sigma$. In this section, we recall some important results concerning the relations between these subspaces and the null-space of the Rosenbrock system matrix pencil $P_{\scriptscriptstyle \Sigma}$.
Given $\mu \in \complex$, we use the symbol $N_{\scriptscriptstyle \Sigma}(\mu)$ to denote a basis matrix for the null-space of $P_{\scriptscriptstyle \Sigma}(\mu)$. 

\subsection{Computation of a basis of $\gR^\star$}

The following result, see \cite{SICON}, presents a procedure for the computation of a basis matrix for $\gR^\star$ and, simultaneously, for the parameterisation of all the friends of $\gR^\star$ that place the eigenvalues of the closed-loop restricted to $\gR^\star$ at arbitrary locations. This procedure aims at constructing a basis for $\gR^\star$ starting from basis matrices $N _{\scriptscriptstyle \Sigma}(\mu_i)$ of the null-spaces of the Rosenbrock matrix relative to the real distinct values $\mu_1,\ldots,\mu_r$, where $r\defi \dim \gR^\star$.

\begin{lemma}{\sc (\cite{SICON})}.
\label{Rs}
Let $r = \dim \gR^\star$. Let $\mu_1,\ldots,\mu_r \in \real$ be distinct and different from the invariant zeros.\footnote{A generalisation of Lemma \ref{Rs} to the case of coincident and possibly complex values of $\mu_i$ is given in \cite{Ntogramatzidis-14}.} 
Let 
\bea
\label{MK}
\bmat{c} \hat{V} \\[-0mm] \hat{W} \emat \defi \bmat{c|c|c|c} N _{\scriptscriptstyle \Sigma}(\mu_1) & N _{\scriptscriptstyle \Sigma}(\mu_2) & \ldots & N _{\scriptscriptstyle \Sigma}(\mu_r) \emat
\eea
be partitioned so that $\hat{V}$ and $\hat{W}$ have $n$ and $m$ rows, respectively.
Then, $\ima \hat{V}=\gR^\star$.
Moreover, let $l$ be the number of columns of $N _{\scriptscriptstyle \Sigma}(\mu_k)$ for any $k \in \{1,\ldots,r\}$.\footnote{The number $l$ does not depend on $k$ since the $\mu_i$ are assumed to be different from the invariant zeros.} Then, 
$V_{\scriptscriptstyle K}\defi \hat{V}\,K$, with $K\defi \diag\{k_1,k_2,\ldots,k_r\}$, is a basis matrix of $\gR^\star$ for almost all $k_i \in \real^{l}$. Finally, defining
$W_{\scriptscriptstyle K}\defi \hat{W}\,K$,
the set of all  friends of $\gR^\star$ such that $\sigma(A+B\,F\,|\,\gR^\star)=\{\mu_1,\ldots,\mu_r\}$ is parameterised in $K$ as 
$F_{\scriptscriptstyle K}=W_{\scriptscriptstyle K}\,V_{\scriptscriptstyle K}^\dagger$, 
where $K$ is such that $\rank V_{\scriptscriptstyle K}=r$.\\
\end{lemma}

Lemma \ref{Rs} permits us to write a spanning set of $\gR^\star$ in terms of the selection of $r$ real numbers.
For any $\mu \in \real \setminus \gZ$, let 
\begin{equation} \label{eq:RsjL}
	\gR^\star (\mu) \stackrel{\text{\tiny def}}{=} \left\{ v \in \gX \,\Big|\,\, \exists w \in \gU \;:\;\;
\left[ \begin{smallmatrix} A-\mu\,I_n & B \\[1mm] C & D \end{smallmatrix} \right]\left[ \begin{smallmatrix} v \\[1mm] w \end{smallmatrix} \right]=0 \right\}.
\end{equation}
From this definition, a decomposition of $\gR^\star$ can be obtained from Lemma~\ref{Rs}. 

\begin{corollary}
\label{agg}
Given distinct $\mu_{1},\cdots,\mu_{r} \in \real \setminus \gZ$, there
holds 
\bea
\label{decomR}
\gR^{\star}=\gR^{\star}(\mu_{1})+\cdots+\gR^{\star}(\mu_{r}).
\eea
\end{corollary}

It is easy to see how to obtain a basis matrix for  $\gR^{\star}(\mu_{i})$. 
Let $\bsmat V_i \\[1mm] W_i \esmat$ be a basis matrix of $\ker P_{\scriptscriptstyle \Sigma}(\mu_i)$ partitioned conformably. Clearly, $\ima V_i=\gR^\star(\mu_i)$. Moreover, $V_i$ is of full column-rank. 
Indeed, if $\omega \in \ker V_i$, then $P_{\scriptscriptstyle \Sigma}(\mu_i)\bsmat V_i \\[1mm] W_i \esmat\,\omega=0$ implies $\bsmat B \\[1mm] D \esmat W_i\,\omega=0$. Since $\bsmat B \\[1mm] D \esmat$ is assumed to be of full column-rank, we conclude that $\omega \in \ker W_i$, so that $\omega\in \ker \bsmat V_i \\[1mm] W_i \esmat=\{0\}$, thus $\omega$ is zero. The matrix $V_i$ is a basis matrix for $\gR^{\star}(\mu_{i})$.

\subsection{Computation of a basis of $\gV^\star_g$}
We now turn our attention to the computation of $\gV^\star_g$. From now on, we will assume that the minimum-phase invariant zeros are all distinct (i.e., their algebraic multiplicity is one):

\begin{assumption}
 \label{Ass2}
System $\Sigma$ has no coincident minimum-phase invariant zeros.
\end{assumption}

This assumption does not lead to a significant loss of generality. In fact, the case of coincident zeros can be dealt with by using the procedure described in \cite{Ntogramatzidis-14}.

The complex numbers $\mu_1,\ldots,\mu_h$  are said to be $s$-conformably indexed if $2\,s \le h$ and $\mu_1,\ldots,\mu_{2\,s}$ are complex, while the remaining are real, and for all odd $k \le 2\,s$ we have $\mu_{k+1}={\mu}_k^\ast$. For example, $\mu_1=1+\i,\mu_2=1-\i,\mu_3=3,\mu_4=-4$ are $1$-conformably indexed; $\mu_1=10\,\i,\mu_2=-10\,\i,\mu_3=2+2\,\i,\mu_4=2-2\,\i$, $\mu_5=7$, are $2$-conformably indexed; $\mu_1=3$, $\mu_2=-1$ are $0$-conformably indexed.
No generality is lost by assuming that, if $\mu_1,\ldots,\mu_h$  are $s$-conformably indexed, for every odd $i \in \{1,\ldots,2\,s\}$, the basis matrix $N _{\scriptscriptstyle \Sigma}(\mu_{i+1})$ is constructed as $N _{\scriptscriptstyle \Sigma}(\mu_{i+1})=N _{\scriptscriptstyle \Sigma}({\mu_{i}}^\ast)={N _{\scriptscriptstyle \Sigma}(\mu_{i})}^\ast$.

\begin{lemma}{\sc (\cite{SICON})}.
\label{cor2}
Let $r = \dim \gR^\star$ and let Assumption~\ref{Ass2} hold. Let $z_{1},\ldots,z_{t}$ be the $s_z$-conformably indexed minimum-phase invariant zeros of $\Sigma$. Let $\mu_1,\ldots,\mu_r$ be $s$-conformably indexed.
 Let also
\beann
{M}_{\scriptscriptstyle \Sigma}(z_k) \defi \left\{ \begin{array}{ll} 
N_{\scriptscriptstyle \Sigma}(z_k)+N_{\scriptscriptstyle \Sigma}(z_{k+1}) & \textrm{if $k\in \{1,\ldots,2\,s_z\}$ is odd} \\[1mm]
i\,[N_{\scriptscriptstyle \Sigma}(z_k)-N_{\scriptscriptstyle \Sigma}(z_{k-1})] & \textrm{if $k\in \{1,\ldots,2\,s_z\}$ is even} \\[1mm]
N_{\scriptscriptstyle \Sigma}(z_k) & \textrm{if $k\in \{2\,s_z+1,\ldots,t\}$}
\end{array} \right.
\eeann
and
\beann
\bmat{c} \hat{V}_g \\ \hat{W}_g \emat \defi  \left[ \begin{array}{c|c|c|c|c|c} \!\! N_{\scriptscriptstyle \Sigma}(\mu_1) & \ldots & N_{\scriptscriptstyle \Sigma}(\mu_r) & M_{\scriptscriptstyle \Sigma}(z_{1}) & \ldots & M_{\scriptscriptstyle \Sigma}(z_{t}) \!\! \end{array} \right].
\eeann
Then, the columns of $\hat{V}_g$ span $\gV^\star_g$. Let $l_k$ be the number of columns of $N_{\scriptscriptstyle \Sigma}(\mu_k)$ and let $\eta_k$ be the number of columns of $M_{\scriptscriptstyle \Sigma}(z_k)$. Then, the columns of the matrix
$V_{\scriptscriptstyle K,H}\defi \hat{V}_g\,\diag\{K,H\}$, with
$K\defi \diag\{k_1,\ldots,k_r\}$ and $H\defi \diag\{h_1,\ldots,h_t\}$ 
are a basis for $\gV^\star_g$ adapted to $\gR^\star$ for almost all $k_i \in \real^{l_i}$ and $h_i \in \real^{\eta_i}$. Finally, defining
$W_{\scriptscriptstyle K,H}\defi \hat{W}_g\,\diag\{K,H\}$,
the set of all  friends $F$ of $\gV^\star_g$ such that the set of eigenvalues of the mapping $A+B\,F$ restricted to $\gV^\star_g$ are $\{\mu_1,\ldots,\mu_r\} \cup \{z_{1},\ldots,z_{t}\}$ 
 is parameterised in $K$ and $H$ as $F_{\scriptscriptstyle K,H}=W_{\scriptscriptstyle K,H}\,V_{\scriptscriptstyle K,H}^\dagger$, 
where $K$ and $H$ are such that $\rank \, V_{\scriptscriptstyle K,H}=\dim \gV^\star_g$.\\[-3mm]
\end{lemma}

\begin{example}
Consider Example \ref{exe0}, for which $r=1$ and which has a minimum-phase invariant zero at $-6$. We choose $\mu_1=-1$, and we compute basis matrices $N_{\scriptscriptstyle \Sigma}(\mu_1)$ and $M_{\scriptscriptstyle \Sigma}(z_1)=N_{\scriptscriptstyle \Sigma}(z_1)$ for $\ker P_{\scriptscriptstyle \Sigma}(\mu_1)$ and $\ker P_{\scriptscriptstyle \Sigma}(z_1)$, respectively:
\beann
\bmat{c} \!\!\! \! \! V_{\scriptscriptstyle K,H}  \! \! \!\!\! \\  \! \! \!\!\! W_{\scriptscriptstyle K,H}   \! \!\!\!\! \emat \!= \! [\begin{array}{cc} 
 \!\!\! N_{\scriptscriptstyle \Sigma}(\mu_1)  \!\!\! &  \!\!\! N_{\scriptscriptstyle \Sigma}(z_1)  \!\!\! \end{array} ]\! \diag\{K,H\} \!=\! \!\bsmat 
 \! 0  \! & \Big| & \! -132 & 0 \\[-3.5mm]
 \! 0  \! & \Big| &  \!44 & 0 \\[-3.5mm]
 \! 0 \!  & \Big| &  \!-123 & 0 \\[-3.5mm]
 \! 1  \! & \Big| &  \!0 & 1 \\[-3.5mm]
 \! 0  \! & \Big| &  \!-6 & 0 \\[-1mm]
\hline \\[-2mm]
 \! -1  \! & \Big| & \! -330 & 0 \\[-3.5mm]
 \! 0  \! & \Big| &  \!216 & 0 \\[-3.5mm]
 \! 0  \! & \Big| & \! 66 & 0 \\[-3.5mm]
 \! 0  \! & \Big| & \! 0 & 0 \esmat \!\diag\{K,H\}
\eeann
where $K$ is $1 \times 1$ and $H$ is $2 \times 1$. For almost all $K$ and $H$ the rank of $V_{\scriptscriptstyle K,H}$ is $2$. Such rank becomes zero only when $K=0$ or the upper entry of $H$ is zero. Choosing e.g. $K=1$ and $H=\bsmat 1 \\[1mm] 0 \esmat$ yields $\gV^\star_g=\ima \bsmat 0 & 0 & 0 & 1 & 0 \\[1mm] -132 & 44 & -123 & 0  & -6 \esmat^\top$.
\end{example}

\section{Solution to {Problem~\ref{prob:mono2}}}

In this section we provide tractable necessary and sufficient conditions for the existence of a solution to the problem of global monotonicity. 
As explained above, in order to achieve a globally monotonic step response we need to find a feedback matrix $F$ that 
renders at least $n-p$ of the $n$ closed-loop modes invisible at the tracking error and evenly distributes the remaining modes evenly into the $p$ components of the tracking error. The number of closed-loop modes that can be made invisible by state feedback equals the dimension of the subspace $\gV^\star_g$. Thus, for the tracking control problem with global monotonicity to be solvable we need the condition $\dim \gV^\star_g \ge n-p$ to be satisfied. This condition is only necessary, because we need also the {linearly independent} vectors $v_1,\ldots,v_p$ obtained with the procedure indicated above to be linearly independent from $\gV^\star_g$. In the case in which $\dim \gV^\star_g > n-p$ holds, if it is possible to find {linearly independent} vectors $v_1,\ldots,v_p$ that are independent from $\gV^\star_g$, then not only is the monotonic tracking control problem solvable, but we are able to also obtain a response that achieves instantaneous tracking in some outputs.

Let $\l \in \real$. For all $j \in \{1,\ldots,p\}$ we define 
\bea
\label{Rhat}
	\hspace{-5mm} \hat{\gR}_j (\lambda) \ns&\ns \stackrel{\text{\tiny def}}{=} \ns&\ns \left\{ v \in \gX \,\Big|\,\, \exists \beta \in \mathbb{R}\setminus \{0\}, \, \exists w \in \gU \;:\;\; \right. \nonumber \\
	 \ns&\ns \ns&\ns \qquad \qquad \quad\left.
\left[ \begin{smallmatrix} A-\lambda\,I_n & B \\[1mm] C & D \end{smallmatrix} \right]\left[ \begin{smallmatrix} v \\[1mm] w \end{smallmatrix} \right]=\left[ \begin{smallmatrix} 0 \\[1mm] \beta e_j \end{smallmatrix} \right] \right\}.
\eea

It is easy to see that, given $\lambda \in \real$, the set $\hat{\gR}_j (\lambda)$ is not a subspace of $\gX$. The set $\hat{\gR}_j (\lambda)$ represents the set of initial states such that a feedback matrix $F$ exists that renders all the output components identically zero with the only exception of the $j$-th component, which must be non-zero. Indeed, given $v \in \hat{\gR}_j (\lambda)$ with $v \neq 0$ and $\beta \in \mathbb{R}\setminus \{0\}$ and $w \in \gU$ such that $\left[ \begin{smallmatrix} A-\lambda\,I_n & B \\[1mm] C & D \end{smallmatrix} \right]\left[ \begin{smallmatrix} v \\[1mm] w \end{smallmatrix} \right]=\left[ \begin{smallmatrix} 0 \\[1mm] \beta e_j \end{smallmatrix} \right]$, the feedback matrices satisfying $F\,v=w$ guarantee that for any $x_0\in \spanR\{v\}$ there hold $y_k=0$ for all $k \in \{1,\ldots,p\} \setminus \{j\}$ and $y_j \neq 0$.

In the following lemma, a first necessary and sufficient condition for the solvability of Problem~\ref{prob:mono2} is given in terms of the existence of vectors $v_{1} \in\hat{\gR}_{1} (\lambda_{1})$,  $\ldots$, $v_{p} \in\hat{\gR}_{p} (\lambda_{p})$ such that $\gV^\star_g + \spanR \{ v_{1},\cdots,v_{p}\} = \gX$. This condition says that 
\begin{enumerate}
\item when  $\dim \gV^\star_g=n-p$, then $v_1,\ldots,v_p$ have to be linearly independent and they all must be independent from $\gV^\star_g$. When this is the case, the real numbers $\lambda_1,\ldots,\lambda_p$ are all part of the closed-loop spectrum for any feedback matrix that solves the problem;
\item when  $\dim \gV^\star_g>n-p$, it may be possible to exploit the excess in ``good'' dimension of $\gV^\star_g$ to compensate for possibly dependent vector(s) $v_i$. In other words, now we do not necessarily need all vectors $v_1,\ldots,v_p$ to be linearly independent and/or independent from $\gV^\star_g$. If this necessary and sufficient condition is satisfied and $\dim \gV^\star_g>n-p$, for any $v_k$ that is dependent on $\gV^\star_g$ or the remaining $v_i$, Problem~\ref{prob:mono2} can be solved with a matrix $F$ such that $\lambda_k$ is not part of the closed-loop spectrum.
\end{enumerate}
\begin{lemma}
\label{lem:Rosej2}
Let Assumptions \ref{Ass1} and \ref{Ass2} hold. Let $\lambda_1, \lambda_2, \cdots, \lambda_p \in \Lambda_g$.
Problem~\ref{prob:mono2} is solvable if and only if for all $j \in \{1,\ldots,p\}$ there exists $v_{j} \in\hat{\gR}_{j} (\lambda_{j})$ such that
\begin{equation} 
\label{eq:Rosej2}
 \gV^\star_g + \spanR \{ v_{1},\cdots,v_{p}\} = \gX.
\end{equation}
\end{lemma}
\proof
Let us consider for the sake of argument the continuous time. The discrete case follows with the obvious substitutions. First, we prove sufficiency. 
Condition (\ref{eq:Rosej2}) guarantees that we can find 
\begin{itemize}
\item $\{v_{\theta_1},\ldots,v_{\theta_l}\}\subseteq \{ v_{1},\cdots,v_{p}\}$, where $\theta : \{1,\ldots,l\} \longrightarrow \{1,\ldots,p\}$ is an injective map; we can define $\hat{v}_i \defi v_{\theta_i}$ for all $i \in \{1,\ldots,l\}$;
\item $\{\tilde{v}_{l+1},\ldots,\tilde{v}_{n}\}\subset \gV^\star_g$,
\end{itemize}
such that $\{v_{\theta_1},\ldots,v_{\theta_l},\tilde{v}_{l+1},\ldots,\tilde{v}_{n}\}=\{\hat{v}_{1},\ldots,\hat{v}_{l},\tilde{v}_{l+1},\ldots,\tilde{v}_{n}\}$ is linearly independent.
Thus, from definition (\ref{Rhat}), there exists $\hat{w}_i\in \gU$ and $\beta_i \neq 0$ such that 
$\bsmat A-\lambda_{\theta_i}\,I_n & B \\[1mm] C & D \esmat \bsmat \hat{v}_i \\[1mm] \hat{w}_i \esmat=
\bsmat 0 \\[1mm] \beta_{\theta_i} \,e_{\theta_i} \esmat$ with $\beta_{\theta_i} \neq 0$ for all $i \in \{1,\ldots,l\}$. Since vectors $\{\tilde{v}_{l+1},\ldots,\tilde{v}_{n}\}$ are such that
$\{\hat{v}_{1},\ldots,\hat{v}_{l},\tilde{v}_{l+1},\ldots,\tilde{v}_{n}\}$ is linearly independent, we can also find vectors $\{\hat{v}_{l+1},\ldots,\hat{v}_n\}$ such that $\{\hat{v}_{1},\ldots,\hat{v}_{l},\hat{v}_{l+1},\ldots,\hat{v}_{n}\}$ is linearly independent by extracting $\{\hat{v}_{l+1},\ldots,\hat{v}_n\}$ from the columns of $V_g = \hat{V}_g\,\diag\{K,H\}$, where $\hat{V}_g$ and the diagonal matrices $K$ and $H$ are constructed as in Lemma \ref{cor2}; in this way, we can define a corresponding set $\{\hat{w}_{l+1},\ldots,\hat{w}_n\}$ from the columns of $W_g = \hat{W}_g\,\diag\{K,H\}$, where $\hat{W}_g$ is also constructed as in Lemma \ref{cor2}.
Let $\nu \defi  \dim \gV^\star_g$. Let us denote by $\delta: \{1,\ldots,\nu\} \longrightarrow \{1,\ldots,n-l\}$ an injective mapping such that $\hat{v}_{l+j}$ is associated with the eigenvalue $\mu_{\delta_j}$, and that $\mu_{\delta_1}, \ldots,\mu_{\delta_{\nu}}$ are $s$-conformably indexed, so that if $\mu_{\delta_j} \in \real$ we have
$\bsmat A-\mu_{\delta_i}\,I_n & B \\[1mm] C & D \esmat \bsmat \hat{v}_{l+i} \\[1mm] \hat{w}_{l+i} \esmat=0$, 
and if $\mathfrak{Re}\{\mu_{\delta_j}\} \neq 0$ with $j$ odd, then $\bsmat \tilde{v}_{l+j+1} \\[1mm] \tilde{w}_{l+j+1} \esmat$ is the complex conjugate of $\bsmat \tilde{v}_{l+j} \\[1mm] \tilde{w}_{l+j} \esmat$. 
Then, the feedback $F \defi [\,\hat{w}_1\,\ldots\,\hat{w}_n\,]\,[\,\hat{v}_1\,\ldots\,\hat{v}_n\,]^{-1}$ satisfies
\beann
&& \hspace{-6mm}(A+B\,F) [\,\hat{v}_1\; \ldots\; \hat{v}_l\;|\;\hat{v}_{l+1}\;\hat{v}_{l+2}\;|\;\ldots\;| \\
&& \qquad \quad \;\hat{v}_{l+2\,s-1}\;\hat{v}_{l+2\,s}\;|\;\hat{v}_{l+2\,s+1} \; \ldots\; \hat{v}_n\,] \\
&&
 =\diag \left\{ \! \lambda_1,\ldots, \lambda_l,\!
 \bsmat \! \mathfrak{Re}\{\mu_1\} \! & -\mathfrak{Im}\{\mu_1\}\! \\[1mm]
 \! \mathfrak{Im}\{\mu_1\}\! & \mathfrak{Re}\{\mu_1\} \! \esmat\!, \ldots,\right. \\
 && \quad \left. \bsmat \! \mathfrak{Re}\{\mu_{2\,s-1}\}\ & -\mathfrak{Im}\{\mu_{2\,s-1}\}\! \\[1mm]
\ \mathfrak{Im}\{\mu_{2\,s-1}\} \! & \mathfrak{Re}\{\mu_{2\,s-1}\} \! \esmat\!, \mu_{2\,s+1},\ldots,\mu_{n-l}\right\}  [\,\hat{v}_1\; \ldots\; \hat{v}_l \\
&& \qquad \quad \;|\;\hat{v}_{l+1}\;\hat{v}_{l+2}\;|\;\ldots\;| \;\hat{v}_{l+2\,s-1}\;\hat{v}_{l+2\,s}\;|\;\hat{v}_{l+2\,s+1} \; \ldots\; \hat{v}_n\,]\\[2mm]
&& \hspace{-4mm} (C+D\,F)\,\hat{v}_i= \left\{ \begin{array}{ll} \beta_{\theta_i}\,e_{\theta_i} & i \in \{1,\ldots,l\} \\
0 & i \in \{l+1,\ldots,n\}
\end{array} \right.
\eeann
 Let $\xi_0=\xi(0)$ be the initial error state, and define $\alpha \defi [\,\hat{v}_1\,\ldots\,\hat{v}_n\,]^{-1}\,\xi_0$. We find
\beann
\epsilon(t) \ns&\ns\! =\! \ns&\ns (C\!+\!D F)\,\exp\left[(A\!+\!B F) t \right] \xi_0 \\
 \ns&\ns\! =\! \ns&\ns \sum_{i=1}^l e^{\lambda_{\theta_i} t}\,(C\!+\!D F) \hat{v}_i \alpha_i= \sum_{i=1}^l \beta_{\theta_i}\,e_{\theta_i}\,\exp(\lambda_{\theta_i}\,t)\,\alpha_i, 
 \eeann
 so that each component of $\epsilon(t)$ is either zero, or it is given by a single exponential in view of the injectivity of $\theta$, and is therefore monotonic.\\
 {
Let us now consider necessity. If Problem~\ref{prob:mono2} admits the solution $F$, by Lemma \ref{lem:mono}, the tracking error has a single closed-loop mode per component, i.e., it is in the form given by (\ref{risp}). This implies that the remaining $n-p$ closed-loop modes (which are asymptotically stable because $F$ is stabilising) must disappear from the tracking error. Hence, $\dim \gV^\star_g \ge n-p$.
Clearly, in (\ref{risp}) some components may be zero due to the fact that the corresponding $\gamma_i$ are zero. Let us assume that the output components are ordered in such a way that the zero components are the last $p-c$, i.e.,
\bea
\epsilon(t) \ns&\ns = \ns&\ns [ \begin{array}{cccccc} \gamma_1\,e^{\lambda_1\,t} & \ldots & \gamma_c\,e^{\lambda_c\,t} & 0  &\ldots & 0 \end{array} ]^\top,
\label{risp1}
\eea
where now $\gamma_i \neq 0$  for $i \in \{1,\ldots,c\}$.
Thus, $\lambda_1,\ldots,\lambda_c\in \sigma(A+B\,F)$. We can partition the matrix $V \defi [\,v_1\;\;v_2\;\;\ldots\;\;v_n\,]$ of the corresponding generalised eigenvectors of $A+B\,F$ as $V=[ \begin{array}{cc} V_1 & V_2 \end{array} ]$, where $V_1=[ \begin{array}{ccc} v_1 & \ldots & v_c \end{array} ]$ has $c$ columns and $\ima V_2\subseteq \gV^\star_g$. On the other hand, we also know that given the initial state $\xi_0 \in \gX$ and $\alpha =V^{-1} \xi_0$ (so that $\xi_0=V_1\,\alpha_1+V_2\,\alpha_2$ when decomposing $\alpha=\bsmat \alpha_1 \\[1mm] \alpha_2 \esmat$ conformably with $V=[ \begin{array}{cc} V_1 & V_2 \end{array} ]$) we can write
\beann
\epsilon(t) \ns&\ns = \ns&\ns (C+D\,F)\,e^{(A+B\,F)\,t}\,\xi_0 \\
\ns&\ns = \ns&\ns (C+D\,F)\,e^{(A+B\,F)\,t}\,(V_1\,\alpha_1+V_2\,\alpha_2).
\eeann
Since $\ima V_2$ is output-nulling for $(A,B,C,D)$ and $F$ is an associated friend,  we have $\ima [e^{(A+B\,F)\,t}\,V_2]\subseteq \ima V_2 \subseteq \gV^\star_g$ and $(C+D\,F)\,e^{(A+B\,F)\,t}\,V_2=0$. Thus
\beann
\epsilon(t) \ns&\ns = \ns&\ns (C+D\,F)\,e^{(A+B\,F)\,t}\,V_1\,\alpha_1 \\
 \ns&\ns = \ns&\ns (C+D\,F)\,V_1\,\diag\{ e^{\lambda_1\,t},\ldots, e^{\lambda_c\,t} \}\,\alpha_1 \\
 \ns&\ns = \ns&\ns (C+D\,F)\,[ \begin{array}{ccc} v_1\,e^{\lambda_1\,t} & \ldots &  v_c\,e^{\lambda_c\,t} \end{array} ]\,\alpha_1.
 \eeann
  As such, for any $i \in \{1,\ldots,c\}$ we have found $w_i \defi F\,v_i$ and $\beta_i \defi \gamma_i \neq 0$ such that $(A-\lambda_i\,I_n)\,v_i +B\,w_i =0$ and $C\,v_i+D\,w_i=\gamma_i\,e_i$ for all $i \in \{1,\ldots,c\}$. This implies that $v_i \in \hat{\gR}_i(\lambda_i)$ for all $i \in \{1,\ldots,c\}$. By taking any vectors $v_{c+1}\in \hat{\gR}_{c+1}(\lambda_{c+1})$, $\ldots$ , $v_p\in \hat{\gR}_{p}(\lambda_{p})$, we have $\gV^\star_g+\spanR\{v_1,\ldots,v_p\}=\gX$.
  }
\endproof
As a direct consequence of the previous result, when Problem~\ref{prob:mono2} is solvable and $\dim  \gV^\star_g=n-p$, all the ``visible'' eigenvalues $\lambda_1,\ldots,\lambda_p$ are included in the closed-loop spectrum for any feedback matrix that solves Problem~\ref{prob:mono2}. This is not necessarily the case when $\dim \gV^\star_g>n-p$, because one or more minimum-phase invariant zeros can be used in place of a vector in one of the sets $\hat{\gR}_{j}(\lambda_j)$, and the corresponding value $\lambda_j$ does not have to be included in the closed-loop spectrum (and the corresponding coefficient $\beta_i$ in (\ref{eq:error}) is equal to zero, which means that in the $i$-th component of the tracking error we achieve exact tracking). \\[-3mm]
\begin{corollary}
Let Assumptions \ref{Ass1} and \ref{Ass2} hold.  If $\dim \gV^\star_g<n-p$, then Problem~\ref{prob:mono2} does not admit solutions.
\end{corollary}
\begin{corollary}
Let Assumptions \ref{Ass1} and \ref{Ass2} hold. 
If Problem~\ref{prob:mono2} is solvable with a feedback matrix $F$ such that 
$\lambda_i \notin \sigma(A+B\,F)$ for some $i \in \{1,\ldots,p\}$, then $\dim \gV^\star_g>n-p$.\\[-3mm]
\end{corollary}

\begin{remark}
Whenever (\ref{eq:Rosej2}) is satisfied, Problem~\ref{prob:mono2} can be solved with an arbitrary convergence rate. At first glance, this property seems to be in contrast with the fact that the pair $(A,B)$ has not been assumed to be completely reachable, but only stabilisable. In other words, one may argue that the uncontrollable modes (which are asymptotically stable), may limit the convergence rate. However, it is easy to see that this is not the case. Indeed, from the right invertibility of the quadruple $(A,B,C,D)$, one can conclude that every uncontrollable
eigenvalue of the pair $(A,B)$ is also an invariant zero of $\Sigma$. \footnote{This can be seen by observing that an uncontrollable eigenvalue $\lambda$ of $(A,B)$ either belongs to
 $\sigma (A+B\,\Phi \,|\,\gX/\gV^\star\!+\!\gR_0)$ or to $\sigma (A+B\,\Phi \,|\,\gV^\star / \gR^\star)$, where $\gR_0=\langle A, \ima B \rangle$ is the reachable subspace of the pair $(A,B)$, i.e., the smallest $A$-invariant subspace containing the range of $B$, and $\Phi$ is any
 friend of $\gV^\star$. Since $\gR_0$ is contained in the smallest input-containing subspace $\gS^\star$ of $\Sigma$ \cite[Chapter 8]{Trentelman-SH-01}, and the right-invertibility is equivalent to the condition $\gV^\star + \gS^\star=\gX$ since  $[\,C\;\;\;D\,]$ has been assumed to be of full row-rank \cite[Theorem 8.27]{Trentelman-SH-01}, we also have $\gV^\star + \gR_0=\gX$. Hence, $\lambda \in \sigma (A+B\,\Phi \,|\,\gV^\star / \gR^\star)$, i.e., $\lambda\in \gZ$. }
Hence, every uncontrollable eigenvalue of the pair $(A,B)$ is rendered invisible at the tracking error, and therefore it does not limit the rate of convergence. 
It is also worth observing that  there is freedom in the choice of the closed-loop eigenvalues associated with $\gR^\star$, when computing a basis matrix for $\gV^\star_g$. Even though these eigenvalues are invisible at the tracking error (and hence any choice will be correct as long as they are asymptotically stable and distinct from the minimum-phase invariant zeros) this freedom may be important for the designer, since the selection of closed-loop eigenvalues affects other considerations like control amplitude/energy. Thus, it is worth emphasising that the designer has complete freedom to chose any set of stable eigenvalues provided the minimum-phase invariant zeros are included, and provided at least $p$ of these meet the desired convergence rate.
\\[-2mm]
\end{remark}

Lemma \ref{lem:Rosej2} already provides a set of necessary and sufficient conditions for the solvability of the globally monotonic tracking control problem. However, such conditions are not easy to test, because they are given in terms of the sets $\hat{\gR}_j(\lambda_j)$ which are not, in general, subspaces of $\gX$.
The tools that we now present are aimed at replacing $\hat{\gR}_j(\lambda_j)$ in condition (\ref{eq:Rosej2}) with particular reachability subspaces of the state-space, which we now define.
{As in the proof of Lemma \ref{lem:Rosej2}, 
for each output} $j\in \{1,\ldots,p\}$ we introduce $\Sigma_j=(A,B,C_{(j)},D_{(j)})$ as the quadruple in which $C_{(j)}\in \real^{(p-1) \times n}$ and $D_{(j)}\in \real^{(p-1) \times m}$ are obtained by eliminating the $j$-th row from $C$ and $D$, respectively.
We observe that the right invertibility of the quadruple $(A,B,C,D)$ guarantees that the set $\gZ$ of invariant zeros of $\Sigma$ contains the set of invariant zeros $\gZ_j$ of $\Sigma_j$ for any $j\in \{1,\ldots,p\}$.
The largest output nulling reachability subspace of $\Sigma_j$ is denoted by $\gR^\star_j$. {Similarly to what was done for $\gR^\star$ in Corollary \ref{agg}, {for any distinct $\mu_{1},\cdots,\mu_{r_j} \in \real \setminus \gZ$,} we decompose $\gR^\star_j$ as}
\begin{equation} \label{eq:gRdec}
\gR^{\star}_j=\gR^{\star}_j(\mu_{1})+\cdots+\gR^{\star}_j(\mu_{r_j}),
\end{equation}
where $r_j = \dim \gR^\star_j$ {and}
\begin{equation} \label{eq:Rsj}
	\gR^\star_j (\mu_i) \stackrel{\text{\tiny def}}{=} \left\{ v \in \gX \,\Big|\,\, \exists w \in \gU \;:\;\;
\left[ \begin{smallmatrix} A-\mu_i\,I_n & B \\[1mm] C_{(j)} & D_{(j)} \end{smallmatrix} \right]\left[ \begin{smallmatrix} v \\[1mm] w \end{smallmatrix} \right]=0 \right\}.
\end{equation}

\begin{remark}
{As established in Corollary \ref{agg}, a spanning set} for $\gR^\star_j (\mu_i)$ is given by the columns of $V_i$, where $V_i$ is the upper part of a basis matrix $\bsmat V_i \\[1mm] W_i \esmat$ of $\left[ \begin{smallmatrix} A-\mu_i\,I_n & B \\[1mm] C_{(j)} & D_{(j)} \end{smallmatrix} \right]$. However, differently from $\gR^\star(\mu)$, this time it is not guaranteed that $V_i$ obtained in this way is of full column-rank, because the matrix $\left[ \begin{smallmatrix} B \\[1mm] D_{(j)} \end{smallmatrix} \right]$, differently from $\left[ \begin{smallmatrix} B \\[1mm] D \end{smallmatrix} \right]$, may very well have a non-trivial kernel. \end{remark}

The relationship between $\hat{\gR}_j (\mu)$ and $\gR^\star_j(\mu)$ is examined in the two following results.

\begin{proposition} 
\label{prop:RhRs}
Let $\mu \in \real \setminus \gZ$. For all $j \in \{ 1,\cdots , p\}$
\begin{equation}
	\gR^\star_j (\mu) = \hat{\gR}_j (\mu) \cup \gR^\star (\mu).
\end{equation}
\end{proposition}

\proof
First, we prove that $\gR^\star_j (\mu) \supseteq \hat{\gR}_j (\mu) \cup \gR^\star (\mu)$. To this end, we first show that $\hat{\gR}_j (\mu) \subseteq \gR^\star_j (\mu)$. Let $v \in \hat{\gR}_j (\mu)$. There exist $w \in \gU$ and $\beta \in \mathbb{R}\setminus \{0\}$ such that
$\left[ \begin{smallmatrix} A-\mu\,I_n & B \\[1mm] C & D \end{smallmatrix} \right]\left[ \begin{smallmatrix} v \\[1mm] w \end{smallmatrix} \right]=\left[ \begin{smallmatrix} 0 \\[1mm] \beta e_j \end{smallmatrix} \right]
$,
which gives $C_{(j)}\,v+D_{(j)}\,w=0$. Hence, $v \in \gR^\star_j (\mu)$.
{We now show that $\gR^\star (\mu) \subseteq \gR^\star_j (\mu)$. Let $v \in \gR^\star (\mu)$. Then, there exist $w \in \gU$ such that
$\left[ \begin{smallmatrix} A-\mu\,I_n & B \\[1mm] C & D \end{smallmatrix} \right]\left[ \begin{smallmatrix} v \\[1mm] w \end{smallmatrix} \right]=\left[ \begin{smallmatrix} 0 \\[1mm] 0 \end{smallmatrix} \right]$,
which again implies that $C_{(j)}\,v+D_{(j)}\,w=0$, so that $v \in \gR^\star_j (\mu)$. Hence, $\hat{\gR}_j (\mu) \cup \gR^\star (\mu) \subseteq \gR^\star_j (\mu)$ holds.}
We show that $\gR^\star_j (\mu) \subseteq \hat{\gR}_j (\mu) \cup \gR^\star (\mu)$. Let
 $v$ be an element of $\gR_j^\star (\mu)$. A $w \in \gU$ exists such that
$\left[ \begin{smallmatrix} A-\mu\,I_n & B \\[1mm] C_{(j)} & D_{(j)} \end{smallmatrix} \right]\left[ \begin{smallmatrix} v \\[1mm] w \end{smallmatrix} \right]=0$. 
Let $\beta = C_{j}\,v+D_{j}\,w$. Then, $\left[ \begin{smallmatrix} A-\mu\,I_n & B \\[1mm] C & D \end{smallmatrix} \right]\left[ \begin{smallmatrix} v \\[1mm] w \end{smallmatrix} \right]=\left[ \begin{smallmatrix} 0 \\[1mm] \beta e_j \end{smallmatrix} \right]$. If $\beta \neq 0$, we have $v \in \hat{\gR}_j (\mu)$, whereas if $\beta=0$, we find $v \in \gR^\star (\mu)$. Thus, $v \in \hat{\gR}_j (\mu) \cup \gR^\star (\mu)$.
\endproof

\begin{proposition} 
\label{prop:meashat}
Let $\mu \in \real \setminus \gZ$. For all $j \in \{ 1,\cdots , p\}$, $\gR^\star_j (\mu)$ and $\hat{\gR}_j(\mu)$ coincide almost everywhere, i.e., they are equal modulo a set defined by the
common zeros of a set of finitely many polynomial equations.
\end{proposition}

\proof
Since $\Sigma$ is right invertible and $\mu$ is not an invariant zero, the inclusion $\gR^\star (\mu) \subseteq \gR_j^\star (\mu)$ deriving from Proposition \ref{prop:RhRs} becomes $\gR^\star (\mu) \subset \gR_j^\star (\mu)$. Indeed, in such a case, $[\, C_{j}\;\;\; D_{j}\,]$ is linearly independent from every row of $\left[ \begin{smallmatrix} A-\mu\,I_n & B \\[1mm] C_{(j)} & D_{(j)} \end{smallmatrix} \right]$. This implies that $\dim \gR^\star (\mu) < \dim \gR_j^\star (\mu)$. Moreover, 
Proposition \ref{prop:RhRs} ensures 
 that $\gR^\star_j (\mu) \setminus \hat{\gR}_j (\mu) \subseteq \gR^\star (\mu)$, which in general does not hold as an equality since $\gR^\star (\mu)$ and $\hat{\gR}_j (\mu)$ may have non-zero intersection. \endproof

 Roughly speaking, this result, {together with Proposition~\ref{prop:RhRs}}, implies that $\gR^\star_j (\mu)$ is coincident with $\hat{\gR}_j (\mu)$ modulo a set of points belonging to a proper algebraic variety within $\gR^\star(\mu)$.  This essential step justifies the fact that from now on we will use 
 $\gR^\star_j (\mu)$, instead of $\hat{\gR}_j (\mu)$, to establish constructive necessary and sufficient condition for our tracking problem.

\begin{example}
Consider Example \ref{exe0}. We $\gR^\star_j(-1)$ ($j =1,2,3$) by computing the basis matrices
\bea
\bmat{c} V_1 \\ W_1 \emat \ns&\ns = \ns&\ns 
\bsmat 0   & 27 &  -80 &  116 &   12  & \Big| &116  &  36 &   36   & 36 \\[-3mm]
 0   &  0 &    0  &   1  &   0  & \Big| &  -1  &   0  &  0   &  0\esmat^\top \label{eqVW1} \\
 \bmat{c} V_2 \\ W_2 \emat\ns&\ns = \ns&\ns 
\bsmat 
 0   &  0  &  28  & -37 &   -6  & \Big| & -37 &  -18 &    0  &   0\\[-3mm]
    0   &  0 &    0  &   1  &   0  & \Big| &  -1  &   0  &  0   &  0\esmat^\top  \label{eqVW2} \\
      \bmat{c} V_3 \\ W_3 \emat\ns&\ns = \ns&\ns 
\bsmat 
 0   &  -27  &  28  & -55 &   -6  & \Big| & -55 &  -18 &    0  &  -36\\[-3mm]
     0   &  0  &   0   &  1  &  0   & \Big| &  -1   & 0  &  0   &  0\esmat^\top  \label{eqVW3} 
      \eea
  of $\ker P_{\scriptscriptstyle \Sigma_j}(-1)$,  which yield $\gR_j^\star(-1)=\ima V_j$.
      \end{example}

\subsection{Solution of Problem~\ref{prob:mono2}: The case $\dim \gV^\star_g=n-p$}

We begin by presenting a famous result in combinatorics \cite[Theorem 3]{Rado-42} due to Rad{\'o}.

\begin{lemma}
Let $P_1,\ldots,P_s$ be sets of a Euclidean space. There exists elements $p_i\in P_i$ for all $i\in \{1,\ldots,s\}$ such that $\{p_1,\ldots,p_s\}$ is a linearly independent set if and only if given $k$ numbers $\nu_1,\ldots,\nu_k$ such that 
$1\le \nu_1 < \nu_2 < \ldots < \nu_k \le s$ for all $k \in \{1,\ldots,s\}$,
the union $P_{\nu_1} \cup P_{\nu_2}\cup \ldots \cup P_{\nu_k}$ contains $k$ linearly independent elements.
\end{lemma} 

\ \\
Let us specialise this result for linear subspaces of $\real^n$. \\

\begin{proposition}
\label{newA}
Let $\gP_1,\ldots,\gP_s$ be subspaces of $\gX$. There exists elements $p_i\in \gP_i$ for all $i\in \{1,\ldots,s\}$ such that $\{p_1,\ldots,p_s\}$ is linearly independent if and only if 
\begin{itemize}
\item $\dim(\gP_{\nu_1})\ge 1$ for all $\nu_1 \in \{1,\ldots,s\}$;
\item $\dim(\gP_{\nu_1}+\gP_{\nu_2})\ge 2$ for all $1\le \nu_1 < \nu_2 \le  s$;\\
\vdots \\[-4mm]
\item $\dim(\gP_{\nu_1}+\gP_{\nu_2}+\ldots+\gP_{\nu_s})\ge s$ for all $1\le \nu_1 < \nu_2 < \ldots < \nu_s  \le  s$.\footnote{This latter condition can be written as $\dim(\gP_1+\ldots+\gP_s)\ge s$.}
\end{itemize}
\end{proposition}
\proof Let $k \in \{1,\ldots,s\}$ and let $1\le \nu_1 < \ldots < \nu_k \le s$. Then, the union $\gP_{\nu_1} \cup \gP_{\nu_2}\cup \ldots \cup \gP_{\nu_k}$ contains $k$ linearly independent elements if and only if the sum $\gP_{\nu_1} + \gP_{\nu_2} + \ldots +  \gP_{\nu_k}$ contains $k$ linearly independent elements, which is in turn equivalent to saying that $\dim (\gP_{\nu_1} + \gP_{\nu_2} + \ldots +  \gP_{\nu_k}) \ge k$. \footnote{Since  $\gP_{\nu_1} +\ldots +  \gP_{\nu_k} \supseteq \gP_{\nu_1} \cup  \ldots \cup \gP_{\nu_k}$,  if $\gP_{\nu_1} \cup  \ldots \cup \gP_{\nu_k}$ contains at least $k$ linearly independent elements, then also $\gP_{\nu_1} +\ldots +  \gP_{\nu_k}$ contains $k$ linearly independent elements. However, the converse is also true because $\gP_{\nu_1} +\ldots +  \gP_{\nu_k}$ is the span of the union $\gP_{\nu_1} \cup  \ldots \cup \gP_{\nu_k}$, so that if $\gP_{\nu_1} \cup  \ldots \cup \gP_{\nu_k}$ had no $k$ linearly independent elements, neither would its span, and therefore there would not exist $k$ linearly independent elements in $\gP_{\nu_1} +\ldots +  \gP_{\nu_k}$.}
\endproof

\begin{corollary}
\label{cornew}
Let $n$ be the dimension of the linear space $\gX$.
Let $\gP_g,\gP_1,\ldots,\gP_s$ be subspaces of $\gX$, and let $\dim \gP_g=n-s$. There exists a linearly independent set $\{p_{g_1},\ldots,p_{g_{n-s}},p_1,\ldots,p_s\}$
such that $\spanR\{p_{g_1},\ldots,p_{g_{n-s}}\}\subseteq \gP_g$ and $p_i\in \gP_i$ for all $i\in \{1,\ldots,s\}$ if and only if 
\begin{itemize}
\item $\dim(\gP_g+\gP_{\nu_1})\ge n-s+1$ for all $\nu_1 \in \{1,\ldots,s\}$;
\item $\dim(\gP_g+\gP_{\nu_1}+\gP_{\nu_2})\ge n-s+2$ for all $1\le \nu_1 < \nu_2 \le  s$; \\
\vdots \\[-4mm]
\item 
$\dim(\gP_g+\gP_1+\ldots+\gP_s)=n$.
\end{itemize}
\end{corollary}
\proof The proof follows directly from the statement of Proposition \ref{newA} with respect to a basis of $\gX$ adapted to $\gP_g$. Indeed, consider $\gX=\gX_1 \oplus \gX_2$ where $\gX_1 = \gP_g$. With respect to this set of coordinates, a basis matrix for $\gP_g$ is given by $\bsmat I_{n-s} \\[1mm] 0_{s \times (n-s)} \esmat$. Let us denote by $\bsmat P_{i,1} \\[1mm] P_{i,2} \esmat$ a basis matrix for $\gP_i$ ($i\in \{1,\ldots,s\}$) with respect to this basis, where $P_{i,1}$ and $P_{i,2}$ have $n-s$ and $s$ rows, respectively. Clearly, we can find   a linearly independent set $\{p_{g_1},\ldots,p_{g_{n-s}},p_1,\ldots,p_s\}$
such that $\spanR\{p_{g_1},\ldots,p_{g_{n-s}}\}\subseteq \gP_g$ and $p_i\in \gP_i$ for all $i\in \{1,\ldots,s\}$ if and only if there exist $\tilde{p}_1 \in \ima P_{1,2}$, $\tilde{p}_2 \in \ima P_{2,2}$, $\ldots$, $\tilde{p}_s \in \ima P_{s,2}$ such that $\{\tilde{p}_1,\ldots,\tilde{p}_s\}$ is linearly independent. However, in view of Proposition \ref{newA} this happens if and only if 
\begin{itemize}
\item $\dim(\gP_{\nu_1,2})\ge 1$ for all $\nu_1 \in \{1,\ldots,s\}$;
\item $\dim(\gP_{\nu_1,2}+\gP_{\nu_2,2})\ge 2$ for all $1\le \nu_1 < \nu_2 \le  s$;\\
\vdots \\[-4mm]
\item $\dim(\gP_{\nu_1,2}+\ldots+\gP_{\nu_s,2})\ge s$ for all $1\le \nu_1 < \ldots < \nu_s  \le  s$.
\end{itemize}
The first condition is equivalent to $\dim(\gP_g+\gP_{\nu_1})\ge n-s+1$ for all $\nu_1 \in \{1,\ldots,s\}$, the second is equivalent to $\dim(\gP_g+\gP_{\nu_1}+\gP_{\nu_2})\ge n-s+2$ for all $1\le \nu_1 < \nu_2 \le  s$, and so on.
\endproof

\ \\
Since in Lemma \ref{lem:Rosej2}
it was shown that when $\dim \gV^\star_g=n-p$  Problem~\ref{prob:mono2} is solvable if and only if there exist 
$v_j \in \hat{\gR}_{j} (\lambda_{j})$ (where $j \in \{1,\ldots,p\}$) satisfying
$\gV^\star_g + \spanR \{ v_{1},\cdots,v_{p}\} = \gX$, and that in Proposition \ref{prop:meashat} it was shown that for any $\mu \in \real\setminus \gZ$ the set $\hat{\gR}_j(\mu)$ coincides with the subspace $\gR^\star_j(\mu)$ modulo a set of points that are roots of an algebraic equation, Corollary \ref{cornew} leads immediately to the following important result.

\begin{theorem} 
\label{th:kchoice}
Let Assumptions \ref{Ass1} and \ref{Ass2} hold.
Let $\dim \gV^\star_g=n-p$.
Let $\lambda_{1},\cdots,\lambda_{p}\in\Lambda_g$.
Problem~\ref{prob:mono2} is solvable if and only if for all ${S}\in2^{\{1,\ldots,p\}}$ there holds
\begin{equation}
\dim\left(\gV_{g}^{\star}+\sum_{j\in{S}}\gR_{j}^{\star}(\lambda_{j})\right)\ge (n-p)+{\rm {card}}({S}). \label{eq:condL}
\end{equation}

\end{theorem}
\ \\
The condition of Theorem \ref{th:kchoice} is a succinct way of writing
$\dim\left(\gV_{g}^{\star}+\gR_{\nu_1}^{\star}(\lambda_{\nu_1})+\ldots+\gR_{\nu_l}^{\star}(\lambda_{\nu_l})    \right)\ge (n-p)+l$ 
for every $1 \le \nu_1  < \ldots < \nu_l \le p$ and every $l \in \{1,\ldots,p\}$, or, explicitly, 
\begin{itemize}
\item $\dim (\gV^\star_g+\gR^\star_j(\lambda_{j}))  \ge  n-p+1$ for all $j \in \{1,\ldots,p\}$; 
\item $\dim (\gV^\star_g+\gR^\star_i(\lambda_{i})+\gR^\star_j(\lambda_{j}))  \ge  n-p+2$ for all
$i,j \in \{1,\ldots,p\}$ such that $i \neq j$;\\
\vdots \\[-4mm]
\item $\dim (\gV^\star_g+\gR^\star_1(\lambda_{1})+\ldots+\gR^\star_p(\lambda_{p}))  =  n$.
\end{itemize}

\begin{example}
Consider the system in Example \ref{exe0}. 
By choosing $\lambda_1=\lambda_2=\lambda_3=-1$ we have $\gR_j^\star(-1)=\ima V_j$, where $V_j$ is given in (\ref{eqVW1}-\ref{eqVW3}).
We recall that $\gV^\star_g=\ima \bsmat 0 && 0 && 0 && 1 && 0 \\[1mm] -132 && 44 && -123 && 0 && -6 \esmat^\top$. Here, (\ref{eq:condL}) is satisfied. Indeed, $\dim(\gV^\star_g+\gR^\star_j(-1))=3$ for all $j \in \{1,2,3\}$, and $\dim (\gV^\star_g+\gR^\star_i(-1)+\gR^\star_j(-1))=4$ for all $i,j \in \{1,2,3\}$ with $i \neq j$. Finally, $\dim(\gV^\star_g+\gR^\star_1(-1)+\gR^\star_2(-1)+\gR^\star_3(-1))=5$.
 \endex
\end{example}

\ \\[-5mm]

We now turn our attention to the problem of the computation of the gain feedback. Assume that the condition in Theorem \ref{th:kchoice} is satisfied, and define $V_{\scriptscriptstyle K,H}\defi \hat{V}_g\,\diag\{K,H\}$ and $W_{\scriptscriptstyle K,H}\defi \hat{W}_g\,\diag\{K,H\}$, where $\hat{V}_g$ and $\hat{W}_g$ are obtained as shown in Lemma \ref{cor2} and $K,H$ are block diagonal matrices constructed as in Lemma \ref{cor2} such that $\ima V_{\scriptscriptstyle K,H}=\gV^\star_g$. 
Matrices $V_{\scriptscriptstyle K,H}$ and $W_{\scriptscriptstyle K,H}$ have $n-p$ columns, and the rank of $V_{\scriptscriptstyle K,H}$ is $n-p$.  Let us partition $V_{\scriptscriptstyle K,H}$ and $W_{\scriptscriptstyle K,H}$ as
$V_{\scriptscriptstyle K,H}=\left[ v_{g,1}\;\;\,v_{g,2}\;\;\ldots \;\;v_{g,n-p}\right]$ and 
$W_{\scriptscriptstyle K,H}=\left[ w_{g,1}\;\;\,w_{g,2}\;\;\ldots \;\;w_{g,n-p}\right]$, which satisfy
$\bsmat A-\mu_i\,I_n & B \\[1mm] C & D \esmat \bsmat v_{g,i} \\[1mm] w_{g,i} \esmat=0$ since $K$ and $H$ are block diagonal (here we assume for the sake of simplicity that all the $\mu_i$ are real).
The necessary and sufficient condition $\gV^\star_g + \spanR \{ v_{1},\cdots,v_{p}\} = \gX$ in Lemma \ref{lem:Rosej2} is satisfied for some {$v_{i} \in \hat{\gR}_{i}(\lambda_{i})$ ($i \in \{1,\ldots,p\}$) because it is assumed that the condition in Theorem \ref{th:kchoice}, which guarantees the solvability of Problem~\ref{prob:mono2}, is satisfied. Hence, there exists $\{w_1,\ldots,w_p\}$ such that $\bsmat A-\lambda_i\,I_n & B \\[1mm] C & D \esmat \bsmat v_{i} \\[1mm] w_{i} \esmat=\bsmat 0 \\ \beta_{i} \,e_i\esmat$ with $\beta_i\neq 0$.} It is now easy to see that the feedback matrix $F_{\scriptscriptstyle K,H}$ that solves Problem~\ref{prob:mono2} is a solution of the linear equation 
\bea
\label{lambda}
F_{\scriptscriptstyle K,H}\,[\begin{array}{ccccc} \! V_{\scriptscriptstyle K,H}  \! &  \! v_1 \!  &  \! \ldots  \! & \! v_p \!  \end{array}] =[\begin{array}{ccccc} \!  W_{\scriptscriptstyle K,H} \!  & \!  w_1 \!  &  \! \ldots \!  & \!  w_p \!  \end{array}].
\eea
Indeed, a feedback matrix $F_{\scriptscriptstyle K,H}$ satisfying (\ref{lambda}) guarantees that 
$(A+B\,F_{\scriptscriptstyle K,H})\,V_{\scriptscriptstyle K,H}=X_g\,V_{\scriptscriptstyle K,H}$ for a certain matrix $X_g$ such that $\sigma(X_g)=\{\mu_1,\ldots,\mu_{n-p}\}$ and $(C+D\,F_{\scriptscriptstyle K,H})\,V_{\scriptscriptstyle K,H}=0$, and that there exist $\beta_i \neq 0$ such that 
$\bsmat A-\lambda_i\,I_n & B \\[1mm] C & D \esmat \bsmat v_{i} \\[1mm] w_{i} \esmat=\bsmat 0 \\ \beta_{i} \,e_i\esmat$, which in turn gives $(A+B\,F_{\scriptscriptstyle K,H})\,v_{i}=\lambda_i\,v_{i}$ and $(C+D\,F_{\scriptscriptstyle K,H})\,v_{i}=\beta_i\,e_i$ for all $i \in \{1,\ldots,p\}$. Therefore, $\sigma(A+B\,F_{\scriptscriptstyle K,H})=\{\lambda_1,\ldots,\lambda_p,\mu_1,\ldots,\mu_{n-p}\}$ and 
\beann
\epsilon(t)\!\! =\!\!(C\!+\!DF_{\scriptscriptstyle K,H}) e^{\lambda_1 t} v_1 \gamma_1\!+\!\ldots\!+\!(C\!+\!DF_{\scriptscriptstyle K,H}) e^{\lambda_p t} v_p \gamma_p \!=\!\!\bmat{c} \!\!\!\!\!\!\!\beta_1\,\gamma_1\,e^{\lambda_1\,t} \!\!\!\!\!\!\!\\[-1mm] \!\!\!\!\!\! \tiny{\vdots}\!\!\!\! \!\!\\ 
\!\!\!\!\!\!\!\beta_p\,\gamma_p\,e^{\lambda_p\,t}\!\!\!\!\!\!\! \emat
\eeann
{for some $\gamma_1,\ldots,\gamma_p \in \real \setminus \{0\}$ as required.}

Since the 
condition $\gV^\star_g + \spanR \{ v_{1},\cdots,v_{p}\} = \gX$ in Lemma \ref{lem:Rosej2}
is equivalent to writing $\rank [\begin{array}{ccccc} V_{\scriptscriptstyle K,H} & v_1 & \ldots & v_p \end{array}]=n$, we can compute a solution $F_{\scriptscriptstyle K,H}$ of (\ref{lambda}) as 
\beann
F_{\scriptscriptstyle K,H} =[\begin{array}{ccccc} \!  W_{\scriptscriptstyle K,H} \!  & \!  w_1 \!  &  \! \ldots \!  & \!  w_p \!  \end{array}]\,\,[\begin{array}{ccccc} \! V_{\scriptscriptstyle K,H}  \! &  \! v_1 \!  &  \! \ldots  \! & \! v_p \!  \end{array}]^{-1}.
\eeann
So far we have shown that if condition (\ref{eq:condL}) in Theorem \ref{th:kchoice} is satisfied, it is possible to find $v_{i} \in \hat{\gR}_{i} (\lambda_{i})$ for all $i \in \{1,\ldots,p\}$ such that $\rank [\,V_g\;\; v_1\;\;\ldots\;\;v_p\,]=n$, and this means that a feedback matrix with the desired properties exists. In other words, so far we have only established the existence of a solution when (\ref{eq:condL}) is satisfied. 
However, a much stronger result holds. Indeed, the vectors $v_1,\ldots,v_p$ can be chosen ``almost randomly'' from within $\gR^\star_1(\lambda_{1})$, $\gR^\star_2(\lambda_{2})$, $\ldots$, $\gR^\star_p(\lambda_{p})$, respectively, and the resulting feedback matrix will almost certainly solve Problem~\ref{prob:mono2} as the following result establishes.\\[-3mm]

\begin{theorem}
\label{theF}
Let Assumptions \ref{Ass1} and \ref{Ass2} hold. Let $\lambda_1,\ldots,\lambda_p\in \Lambda_g$.  Let $r=\dim \gR^\star$ and  $\dim \gV^\star_g =n-p$.
Assume the condition in Theorem \ref{th:kchoice} holds. 
Let $\hat{V}_g$ and $\hat{W}_g$ be constructed as in Lemma \ref{cor2} for the asymptotically stable complex numbers $\mu_1,\ldots,\mu_r$ and for the minimum phase invariant zeros $z_1,\ldots,z_t$.
Let $\bsmat V_i \\[1mm] W_i \esmat$ denote a basis matrix for the kernel of $P_{\scriptscriptstyle \Sigma_i}(\lambda_i)=\bsmat A-\lambda_i\,I_n && B \\[1mm] C_{(i)} && D_{(i)} \esmat$ for all $i \in \{1,\ldots,p\}$, where each
$V_i$ and $W_i$ have $n$ and $m$ rows, respectively. Finally, let
\beann
V_{\scriptscriptstyle K,H,k_1,\ldots,k_p} \ns&\ns \defi \ns&\ns [\begin{array}{ccccc} \hat{V}_g\,\diag\{K,H\} & V_1\,k_1 &V_2\,k_2 & \ldots & V_p\,k_p \end{array} ], \\
W_{\scriptscriptstyle K,H,k_1,\ldots,k_p} \ns&\ns \defi \ns&\ns [\begin{array}{ccccc} \hat{W}_g\,\diag\{K,H\} & W_1\,k_1 &W_2\,k_2 & \ldots & W_p\,k_p \end{array} ],
\eeann
where $k_1,\ldots,k_p \neq 0$ are real parameter vectors of appropriate sizes and $K$ and $H$ are block diagonal parameter matrices as in Lemma \ref{cor2} such that $\ima (\hat{V}_g\,\diag\{K,H\})=\gV^\star_g$. Then:
\begin{itemize}
\item The rank of $V_{\scriptscriptstyle K,H,k_1,\ldots,k_p}$ is equal to $n$ for almost all $K$ and $H$ and $k_1, \ldots,k_p \neq 0$ as constructed above;
\item The set of all feedback matrices that solve Problem~\ref{prob:mono2} for the given $\mu_1,\ldots,\mu_r$ is given by
\bea
\label{param0}
F_{\scriptscriptstyle K,H,k_1,\ldots,k_p}=W_{\scriptscriptstyle K,H,k_1,\ldots,k_p}\,V_{\scriptscriptstyle K,H,k_1,\ldots,k_p}^{-1}.
\eea
\end{itemize}
\end{theorem}
\proof 
In view of Theorem \ref{th:kchoice} and Lemma~\ref{lem:Rosej2}, there exist
$v_{j} \in\hat{\gR}_{j} (\lambda_{j})$ such that \eqref{eq:Rosej2} holds. Since $\hat{\gR}_{j} (\lambda_{j})\subseteq \gR_{j}^\star (\lambda_{j})$ holds by Proposition~\ref{prop:RhRs}, there exist real vectors $k_1,\ldots,k_p$ of suitable sizes such that $\rank \,\mathfrak{A} = n$, where
\bea
\label{appoggio}
\mathfrak{A} \defi
[ \begin{array}{c|c|c|c}
\! \! \hat{V}_g\diag\{K,H\} & \! V_1 & \! \ldots &\! V_p\! \! \end{array} ]
\diag\{I_{n-p},k_1,\ldots,k_p\}.
\eea
Since $\rank [ \begin{array}{c|c|c|c}
\hat{V}_g\,\diag\{K,H\} & V_1 & \ldots & V_p \end{array} ]=n$ in view of (\ref{eq:condL}) written for $S=\{1,\ldots,p\}$, we conclude that $\mathfrak{A}$ loses rank only for $k_1,\ldots,k_p$ that solve a finite set of linear equations.

It remains to show that the parameterisation (\ref{param0}) of the feedback matrices  which solve Problem~\ref{prob:mono2} is exhaustive, i.e., that given a feedback $F$ which solves Problem~\ref{prob:mono2} for $\lambda_1,\ldots,\lambda_p$, there exist $H$, $K$, $k_1$, $\ldots$, $k_p$ such that, computing $V_{\scriptscriptstyle K,H,k_1,\ldots,k_p}$ and $W_{\scriptscriptstyle K,H,k_1,\ldots,k_p}$ as in the statement, $F$ can be written as $W_{\scriptscriptstyle K,H,k_1,\ldots,k_p}\,V_{\scriptscriptstyle K,H,k_1,\ldots,k_p}^{-1}$. In view of Lemma \ref{lem:Rosej2}, $F$ satisfies 
\bea
\label{eq311}
\bmat{c} A+B\,F \\ C+D\,F \emat V_g =\bmat{c} V_g \\[0mm] 0 \emat\,X_g,
\eea
where $V_g$ is a basis matrix for $\gV^\star_g$ and where $X_g$ is asymptotically stable, and 
 \bea
 \label{eq321}
 \bmat{c} A+B\,F \\[0mm] C_{(i)}+D_{(i)}\,F \emat v_i =\bmat{c} v_i \\[0mm] 0 \emat\,\lambda_i,
 \eea
  with $i \in \{1,\ldots,p\}$, where $v_i \in \hat{\gR}_i(\lambda_i)$. Assuming for simplicity that all the eigenvalues of $X_g$ are real and distinct,\footnote{The case of complex eigenvalues of $X_g$ can be dealt with using the argument in \cite[Theorem 3.1]{SICON}.} we can find a change of coordinate matrix $T$ in $\gX$ such that $T^{-1}\,X_g\,T=X_{\scriptscriptstyle \triangle}$ is diagonal. Thus, denoting by $\upsilon_i$ the $i$-th column of $V_g\,T$, and by $\{\mu_1,\ldots,\mu_{n-p}\}$ the eigenvalues of $X_{\scriptscriptstyle \triangle}$, (\ref{eq311}) yields
 \bea
 && \hspace{-1cm} \bmat{c} A+B\,F \\ C+D\,F \emat [\begin{array}{cccc}  \! \! \upsilon_1 \! \! & \! \! \upsilon_2 \! \! & \! \! \ldots \! \! & \! \! \upsilon_{n-p} \! \! \end{array}] \nonumber \\
 && \hspace{-0.4cm}=\bmat{cccc} \! \!  \upsilon_1 \! \! & \! \! \upsilon_2 \! \! & \! \! \ldots \! \! & \! \! \upsilon_{n-p}\! \!   \\[0mm]\! \!  0 \! \! & \! \! 0 \! \! & \! \! \ldots \! \! & \! \! 0\! \! \emat\,\diag \{\mu_{\eta(1)},\ldots,\mu_{\eta(n-p)}\},\label{eq312}
\eea
where $\eta:\{1,\ldots,n\!-\!p\}\! \longrightarrow\! \{1,\ldots,n\!-\!p\}$ is a bijection.
Defining $\omega_i \defi F\,\upsilon_i$, we find that $\bsmat \upsilon_i \\[1mm] \omega_i \esmat \in \ker \bsmat A-\mu_{\eta_i}\,I_n & B \\[0mm] C & D \esmat$. We can repeat the same argument for (\ref{eq321}) (without the diagonalisation), and defining $w_i=F\,v_i$, there holds 
 $\bsmat v_i \\[1mm] w_i \esmat \in \ker \bsmat A-\lambda_i\,I_n & B \\[0mm] C_{(i)} & D_{(i)} \esmat$. Thus, {\bf (i)} $F$ satisfies $[\begin{array}{ccccccc} \!\! \!\omega_1 \! \!\!& \!\! \!\ldots\!\! \! & \!\! \!\omega_{n-p}\!\! \! & \!\! \!w_1 \!\! \!&\! \!\! \ldots\! \! \!&\!\! \! w_p\! \! \!\end{array}]=F\,[\begin{array}{ccccccc} \! \!\! \upsilon_1 \!\! \! & \!\! \! \ldots\! \! \! & \!\! \! \upsilon_{n-p}\! \!\! & \!\! \!v_1 \! \!\!&\! \!\! \ldots \!\! \!&\!\! \! v_p \! \! \!\end{array}]$; {\bf (ii)} $[\begin{array}{ccccccc} \! \!\omega_1 \! \!& \! \!\ldots\! \! & \! \!\omega_{n-p}\! \! & \! \!w_1 \! \!&\! \! \ldots \! \!&\! \! w_p \! \! \end{array}]$ can be written as $W_{\scriptscriptstyle K,H,k_1,\ldots,k_p}$ for a suitable choice of the parameter matrices $K$, $H$ and $k_i$; {\bf (iii)} $[\begin{array}{ccccccc} \! \! \upsilon_1 \! \! & \! \! \upsilon_2 \! \! & \! \! \ldots \! \! & \! \! \upsilon_{n-p}\! \! & \! \!v_1 \! \!&\! \! \ldots \! \!&\! \! v_p \! \! \end{array}]$ can be written as $V_{\scriptscriptstyle K,H,k_1,\ldots,k_p}$ for suitable values of $K$, $H$ and $k_i$. Thus, $W_{\scriptscriptstyle K,H,k_1,\ldots,k_p}=F\,V_{\scriptscriptstyle K,H,k_1,\ldots,k_p}$. \\
\endproof

\begin{example}
Consider again the system in Example \ref{exe0}. Choosing $K=1$, $H=k_1=k_2=k_3=\bsmat 1 \\[1mm] 0 \esmat$, we obtain
\beann
V_{\scriptscriptstyle K,H,k_1,\ldots,k_p}\!=\!\!\! \bsmat 
\!\! 0 \! &  -132  \! &    0 \! &     0  \! &    0\!\! \\[1mm]
 \!\!     0  \! &   44 \! &    27 \! &     0 \! &   -27\!\! \\[1mm]
 \!\!     0 \! &  -123 \! &   -80  \! &   28  \! &   28\!\! \\[1mm]
  \!\!    1  \! &    0  \! &  116 \! &   -37  \! &  -55\!\! \\[1mm]
   \!\!   0  \! &   -6   \! &  12  \! &   -6   \! &  -6\!\! \esmat\!\!\!,
 W_{\scriptscriptstyle K,H,k_1,\ldots,k_p}\!=\!\!\! \bsmat
\!\!  -1\! &   -330  \! &  116 \! &   -37  \! &  -55\!\! \\[1mm]
 \!\!     0 \! &   216  \! &   36 \! &   -18  \! &  -18\!\! \\[1mm]
  \!\!    0  \! &   66 \! &    36 \! &     0   \! &   0\!\! \\[1mm]
  \!\!    0  \! &    0  \! &   36  \! &    0 \! &   -36\!\! \esmat
     \eeann
     We can compute the feedback using (\ref{param0}), which yields (\ref{f1}). The set of parameters $K,H,k_1,k_2,k_3$ for which $V_{\scriptscriptstyle K,H,k_1,\ldots,k_p}$ is singular is given by $K=0$, or $H=\bsmat 0 \\[1mm] \star \esmat$, or $k_1=\bsmat 0 \\[1mm] \star \esmat$, or $k_2=\bsmat 0 \\[1mm] \star \esmat$, or $k_3=\bsmat 0 \\[1mm] \star \esmat$, which constitute a $4$-dimensional algebraic variety in the parameter space which is $9$ dimensional.
        \end{example}

\begin{example}
\label{esprimo}
Consider the right invertible quadruple $(A,B,C,D)$  given by
\beann
A = \bsmat 0 && 0 && 0 && -1  \\[1mm] 0 && -2 && -2 && 0  \\[1mm] 1 && 0 && -4 && 0  \\[1mm]  0 && 0 && 0 && 1 \esmat, \quad
B = \bsmat  2 &&    2&&     0\\[1mm]
     0  &&   0   &&  0\\[1mm]
    -2  &&  0   &&  0\\[1mm]
     0  &&   0   &&  4 \esmat, \quad
C = \bsmat -1 && 0 && 1 && 0 \\[1mm] 0 && -1 && 2 && 0 \esmat
\eeann
and $D=0$. 
The only invariant zero is $z=-3$. The null-space of $P_{\scriptscriptstyle \Sigma}(-3)$ is given by $\ker P_{\scriptscriptstyle \Sigma}(-3)=\ima \bsmat \hat{V}_g \\[1mm] \hat{W}_g \esmat$ where 
\beann
\hat{V}_g=  \bsmat 3 && 0 \\[1mm] 6 && 0 \\[1mm]  3 && 0\\[1mm] 1 && -2\esmat, \quad
\hat{W}_g=  \bsmat  0 && 0 \\[1mm] -4 && -1\\[1mm] -2 && 2 \esmat,
\eeann
(so that $\dim \gV^\star_g=n-p=2$) and taking e.g.  $H=\bsmat 1/3 && 0 \\[1mm] 1/6 && -1/2 \esmat$ and $K$ to be the empty matrix (because in this case $\gR^\star=\{0\}$) leads to 
\beann
\hat{V}_g \,H=\bmat{cc}      
\!\!\!\! 1   \!\! & \!\!   0   \!\!\!\! \\        
\!\!\!\!  2   \!\! & \!\!   0   \!\!\!\! \\        
\!\!\!\!  1   \!\! & \!\!   0   \!\!\!\! \\       
\!\!\!\!   0   \!\! & \!\!   1      \!\!\!\! 
      \emat, \qquad
     \hat{W}_g \,H=\bmat{cc}      
\!\!\!\! 0   \!\! & \!\!   0   \!\!\!\! \\        
\!\!\!\!  -\frac{3}{2}   \!\! & \!\!  \frac{1}{2}   \!\!\!\! \\        
\!\!\!\!   0   \!\! & \!\!   -1      \!\!\!\! 
      \emat.
      \eeann
Consider $\lambda_1=\lambda_2=-1$. We find that basis matrices for $\ker P_{\scriptscriptstyle \Sigma_1}(\lambda_1)$ and $\ker P_{\scriptscriptstyle \Sigma_2}(\lambda_2)$ are respectively given by
\beann
\bmat{c} V_1 \\ W_1 \emat=\bsmat 3 \\[1mm] 0 \\[1mm] 0 \\[1mm] 1 \\[1mm] \hline 3/2 \\[1mm] -5/2 \\[1mm] -1/2 \esmat   \quad \text{and} \quad  \bmat{c} V_2 \\ W_2 \emat=\bsmat -6 \\[1mm] 12 \\[1mm] -6 \\[1mm] 1 \\[1mm] \hline 6 \\[1mm] -5/2 \\[1mm] -1/2 \esmat.
\eeann
Thus, $\gR_1^\star(\lambda_1)=\spanR \{\bmat{cccc} 3 & 0 & 0 & 0 \emat^\top \}$ and $\gR_2^\star(\lambda_2)=\spanR \{\bmat{cccc} -6& 12& -6& 1 \emat^\top \}$. Condition (\ref{eq:condL}) is fulfilled since $\dim(\gV^\star_g+\gR_1^\star(\lambda_1))=\dim(\gV^\star_g+\gR_2^\star(\lambda_2))=3$ and $\dim(\gV^\star_g+\gR_1^\star(\lambda_1)+\gR_2^\star(\lambda_2))=4$. Taking $k_1=1$ and $k_2=2$ gives
\beann
F_{\scriptscriptstyle H,k_1,k_2} \ns&\ns = \ns&\ns [ \begin{array}{ccc} \hat{W}_g \,H & W_1\,k_1 & W_2\,k_2 \end{array} ]\, [ \begin{array}{ccc} \hat{V}_g \,H & V_1\,k_1 & V_2\,k_2 \end{array} ]^{-1} \\
\ns&\ns = \ns&\ns \bsmat 1/2 && 1/4 && -1 && 0 \\[1mm] -1 && -1/2 && 1/2 && 1/2 \\[1mm] 1/6 && -1/16 && -3/8 && -1 \esmat.
\eeann
It is straightforward to check that this feedback matrix solves Problem~\ref{prob:mono2}. 
\endex
\end{example}

\subsection{Solution of Problem~\ref{prob:mono2}: the general case}

We now consider the case where the dimension of $\gV^\star_g$, which we denote by $h$, is possibly strictly greater than  $n-p$. The following generalisation of Rad{\'o}'s Theorem, see \cite[Theorem 1.3]{Mirsky-P-67}, is the key to obtaining a necessary and sufficient solvability condition Problem~\ref{prob:mono2} in this general case.

\begin{proposition}
\label{newB}
Let $\gP_1,\ldots,\gP_s$ be subspaces of $\gX$. There exists $k$ elements $p_{1}\in \gP_{i_1}$, $p_{2}\in \gP_{i_2}$, $\ldots$, $p_{k}\in \gP_{i_k}$ for some $1 \le i_1 < \ldots < i_k \le s$ such that $\{p_1,\ldots,p_k\}$ is linearly independent if and only if 
\begin{itemize}
\item $\dim(\gP_{\nu_1}+\ldots+\gP_{\nu_{s-k+1}})\ge n-k$ for all $1 \le \nu_1< \ldots < \nu_{s-k+1}\le s$;\\
\vdots \\[-4mm]
\item 
$\dim(\gP_{1}+\ldots+\gP_{s})= k$.
\end{itemize}
\end{proposition}

As a result of Proposition \ref{newB}, following the same argument of the proof of Corollary \ref{cornew}, one easily sees that a necessary and sufficient condition for Problem~\ref{prob:mono2} is given in the general case $\dim \gV^\star_g \ge n-p$ as follows.

\begin{theorem}
\label{th:kchoicegen}
Let Assumptions \ref{Ass1} and \ref{Ass2} hold.
Let $\lambda_{1},\cdots,\lambda_{p}\in\Lambda_g$.
Problem~\ref{prob:mono2} is solvable if and only if 
\bea
&&\hspace{-11mm}\dim\left(\gV_{g}^{\star}+\sum_{j\in{S}}\gR_{j}^{\star}(\lambda_{j})\right)\ge n-p+{\rm {card}({S})} \nonumber \\
&& \forall\,{S}\in \{\mathfrak{S} \in 2^{\{1,\ldots,p\}}\,|\;\; \mathrm{card}\,{\mathfrak{S}}>h-(n-p)\}.\label{eq:condLhglobgen}
\eea
\end{theorem}

It is clear that (\ref{eq:condLhglobgen}) reduces to (\ref{eq:condL}) when $h=\dim \gV^\star_g=n-p$. Observe also that (\ref{eq:condLhglobgen}) can be alternatively written as 
\bea
\dim(\gV^\star_g+\gR^\star_{\nu_1}(\lambda_1)+\ldots+\gR^\star_{\nu_l}(\lambda_{\nu_l}))\ge (n-p)+l \label{eq:condLbisgen}
\eea
for every $1 \le \nu_1 < \nu_2 < \ldots < \nu_l \le p$ and every $l \in \{h-(n-p),\ldots,p\}$. 

The calculation of the feedback matrix does not change significantly with respect to the one outlined in Theorem~\ref{theF} for the case $\dim \gV^\star_g = n-p$. The main difference is that 
the $n \times (h+p)$ matrix 
$\mathfrak{V} \defi [\begin{array}{ccccc} \hat{V}_g \diag\{K,H\} & V_1\,k_1 & \ldots & V_p \,k_p \end{array}]$
is not full column-rank. On the other hand, the rank of $\mathfrak{V}$ is $n$ for suitable values of the parameter matrices, which means that it is sufficient to eliminate from $\mathfrak{V}$ exactly $h+p-n$ columns that are linearly
 dependent upon the remaining $n$ columns.

We eliminate the corresponding columns of $[\begin{array}{ccccc}\!\!\!\! \hat{V}_g\diag\{K,H\}\!\!\!\!  & \!\!V_1\!\!\!\! & \!\!\ldots \!\!\!\! & V_p \!\!\!\! \end{array}]$ and $[\begin{array}{ccccc} \!\!\!\! \hat{W}_g\diag\{K,H\}\!\!\!\!  &\!\! W_1\!\!\!\!  &\!\! \ldots \!\!\!\! & W_p\!\!\!\!  \end{array}]$, and we also eliminate the corresponding columns and
 rows from the parameter matrix $\diag\{I_h,k_1,k_2,\ldots,k_p\}$. We denote the matrices thus obtained by $\tilde{V}_{\scriptscriptstyle K,H,k_1,\ldots,k_p,\psi}$, $\tilde{W}_{\scriptscriptstyle K,H,k_1,\ldots,k_p,\psi}$ and $\tilde{K}_{\scriptscriptstyle
 K,H,k_1,\ldots,k_p,\psi}$ respectively, where $\psi$ is a mapping that represents the choice of the columns that have been eliminated.\footnote{For example, if we choose to eliminate the last $h+p-n$ columns of $\mathfrak{V}$, we get $\tilde{V}_{\scriptscriptstyle
 K,H,k_1,\ldots,k_p,\psi} = \begin{bmatrix} \hat{V}_g\diag\{K,H\} & V_1& \ldots & V_{n-h} \end{bmatrix}$, $\tilde{W}_{\scriptscriptstyle K,H,k_1,\ldots,k_p,\psi} = \begin{bmatrix} \hat{W}_g\diag\{K,H\} & W_1& \ldots & W_{n-h} \end{bmatrix}$ and $\tilde{K}_{\scriptscriptstyle
 K,H,k_1,\ldots,k_p,\psi}=\diag\{I_h,k_1,\ldots,k_{n-h}\}$.}
Now, the argument of Theorem \ref{theF} can be applied to the equation 
\bea
\label{paramreduced}
F_{\scriptscriptstyle K,H,k_1,\ldots,k_p,\psi}\,\tilde{V}_{\scriptscriptstyle K,H,k_1,\ldots,k_p,\psi}=\tilde{W}_{\scriptscriptstyle K,H,k_1,\ldots,k_p,\psi},
\eea
which gives the solution to Problem~\ref{prob:mono2} in parameterised form. We have just proved the following result.\\[-3mm]

\begin{theorem}
\label{theF1} 
Let Assumptions \ref{Ass1} and \ref{Ass2} hold. Let $\lambda_1,\ldots,\lambda_p\in \Lambda_g$.  Let $r=\dim \gR^\star$, and let $h = \dim \gV^\star_g  \ge n-p$.
Let the condition in Theorem \ref{th:kchoicegen} hold. 
Let $\hat{V}_g$ and $\hat{W}_g$ be constructed as in Lemma \ref{cor2} for the asymptotically stable complex numbers $\mu_1,\ldots,\mu_r$ and for the minimum phase invariant zeros $z_1,\ldots,z_t$.
Let $\bsmat V_i \\[1mm] W_i \esmat$ denote a basis matrix for the kernel of $P_{\scriptscriptstyle \Sigma_i}(\lambda_i)=\bsmat A-\lambda_i\,I_n && B \\[1mm] C_{(i)} && D_{(i)} \esmat$ for all $i \in \{1,\ldots,p\}$, where each
$V_i$ and $W_i$ have $n$ and $m$ rows, respectively. Finally, let $\nu_i$ denote the number of columns of $V_i$. Let $k_i \in \real^{\nu_i}$ denote a parameter vector for all $i \in \{1,\ldots,p\}$. 
Let $\psi_1,\ldots,\psi_n$ be indexes of the columns of $V = [\begin{array}{ccccc} \hat{V}_g \diag\{K,H\} & V_1& \ldots & V_p  \end{array}]$ such that the rank of 
$\tilde{V}_{\scriptscriptstyle K,H,k_1,\ldots,k_p,\psi}=[\begin{array}{cccc} V^{\scriptscriptstyle \psi_1} & V^{\scriptscriptstyle \psi_2} & \ldots & V^{\scriptscriptstyle \psi_n} \end{array}]$ is equal to $n$. \footnote{We recall that we denote by $V^{\scriptscriptstyle \psi_i}$ the $\psi_i$-th column of $V$.} Let $\tilde{W}_{\scriptscriptstyle K,H,k_1,\ldots,k_p,\psi}=[\begin{array}{cccc} W^{\scriptscriptstyle \psi_1} & W^{\scriptscriptstyle \psi_2} & \ldots & W^{\scriptscriptstyle \psi_n} \end{array}]$, and let $\tilde{K}_{\scriptscriptstyle K,H,k_1,\ldots,k_p,\psi}$ be obtained by $\diag\{I_h,k_1,\ldots,k_p\}$ by removing the corresponding rows and columns. Then:
\begin{itemize}
\item the rank of $\tilde{V}_{\scriptscriptstyle K,H,k_1,\ldots,k_p,\psi}$ is equal to $n$ for almost all $K$ and $H$ as defined in Lemma \ref{cor2}, for all $k_1, \ldots,k_p \neq 0$, and for all the choices $\psi$ such that the matrix obtained by eliminating $h+p-n$ columns from $V$ gives a matrix of rank $n$;
\item The set of feedback matrices that solve Problem~\ref{prob:mono2} with $\mu_1,\ldots,\mu_r$, $z_1,\ldots,z_t$ and the given choice of $\psi$ is 
\bea
\label{param}
F_{\scriptscriptstyle K,H,k_1,\ldots,k_p,\psi}=\tilde{W}_{\scriptscriptstyle K,H,k_1,\ldots,k_p,\psi}\,\tilde{V}_{\scriptscriptstyle K,H,k_1,\ldots,k_p,\psi}^{-1}
\eea
where $k_i \in \real^{\nu_i}$ ($i \in \{1,\ldots,p\}$), $K$ and $H$ are block-diagonal matrices constructed as in Lemma \ref{cor2} such that $\ima (\hat{V}_g\,\diag\{K,H\})=\gV^\star_g$, and $\psi$ is such that the matrix obtained by eliminating $h+p-n$ columns from $V$ gives a matrix of rank $n$. 
\end{itemize}
\end{theorem}

\begin{remark}
If $h=\dim \gV^\star_g > n-p$, the most natural choice is to build a feedback using a basis of the state space that uses as many basis vectors as possible from $\gV^\star_g$, because every extra basis vector (in addition to the first $n-p$) that we use of $\gV^\star_g$ potentially
 results in a tracking error with a further identically zero component. This corresponds to selecting a mapping $\psi$ which eliminates as many columns $V_i \,k_i$  from $\mathfrak{V}$ as possible (under the constraint $\rank \,\tilde{V}_{\scriptscriptstyle K,H,k_1,\ldots,k_p,\psi}
 = n$).
\end{remark}

\begin{example}
\label{esaggiunto2}
Consider the quadruple 
$(A,B,C,D)$ given by
\beann
A = \bsmat -10 & 0 & 8  \\[1mm] 0 & -9 & 0  \\[1mm] 1 & 3 & 10  \esmat\!\!, \quad
B = \bsmat  0&     0\\[1mm]
     0  &   0 \\[1mm] 
  -8 &   0 \esmat\!\!, \quad
C = \bsmat -3 & 0 & 8\\[1mm]  9 & -5 & 6 \esmat\!\!, \quad
D = \bsmat 0 & 0 \\[1mm]  0 & 4 \esmat\!\!,
\eeann
which has two invariant zeros $z_1=-7$ and $z_2=-9$ and is right and left invertible (so that $\gR^\star=\{0\}$). We find
\beann
\ker P_{\scriptscriptstyle \Sigma}(-7)\ns&\ns = \ns&\ns \ima [\begin{array}{ccc|cc}  8&  0 & 3 &  \frac{59}{8}& -\frac{45}{2}  \end{array}]^\top,\\
\ker P_{\scriptscriptstyle \Sigma}(-9)\ns&\ns = \ns&\ns\ima  [\begin{array}{ccc|cc} 0& 1 & 0& \frac{3}{8}&\frac{5}{4}   \end{array}]^\top,
\eeann
so that we can take $V_g=\bsmat 8 && 0 \\[1mm] 0 && 1 \\[1mm] 3 && 0 \esmat$ and $W_g=\bsmat  \frac{59}{8}   &   \frac{3}{8}  \\[1mm]-\frac{45}{2}   &  \frac{5}{4} \esmat$.  Let $\lambda_1=\lambda_2=-2$. We compute basis matrices for $\ker P_{\scriptscriptstyle \Sigma_1}(\lambda_1)$ and $\ker P_{\scriptscriptstyle \Sigma_2}(\lambda_2)$, which are respectively given by
\beann
\bsmat V_1 \\[1mm] W_1 \esmat\ns&\ns = \ns&\ns \ima [\begin{array}{ccc|cc}  -8 &0& -8 & -13 &30 \end{array}]^\top\\ 
\bsmat V_2 \\[1mm] W_2 \esmat\ns&\ns = \ns&\ns  [\begin{array}{ccc|cc}0&0& 0&0& 1 \end{array}]^\top.
\eeann
Thus, $\gR_1^\star(\lambda_1)=\spanR \{\bmat{cccc} 1 & 0 & 1 \emat^\top \}$ and $\gR_2^\star(\lambda_2)=\{0\}$. Condition (\ref{eq:condLhglobgen}) is fulfilled since $\dim(\gV^\star_g+\gR_1^\star(\lambda_1)+\gR_2^\star(\lambda_2))=3$. Taking $H=I_2$, $k_1=1$ and $k_2=0$ ($K$ is the empty matrix since $\gR^\star=\{0\}$) gives
\beann
F_{\scriptscriptstyle K,H,k_1,k_2,\psi} \ns&\ns = \ns&\ns [ \begin{array}{ccc} \hat{W}_g \,H & W_1\,k_1\end{array} ]\, [ \begin{array}{ccc} \hat{V}_g \,H & V_1\,k_1 \end{array} ]^{-1} \\
\ns&\ns = \ns&\ns \bsmat \frac{59}{8} && \frac{3}{8} && -13  \\[1mm] -\frac{45}{2} && \frac{5}{4} && 30   \esmat\!\! \bsmat 8 && 0 && -8  \\[1mm] 0 && 1 && 0  \\[1mm] 3 && 0 && -8  \esmat^{-1}\!= 
 \frac{1}{8}  \bsmat 4 && 3 &&  9 \\[1mm] -18 && 10 && -12 \esmat\!,
\eeann
where $\psi$ is associated with the particular choice $k_2=0$.
\end{example}

\begin{example}
\label{esaggiunto3}
Consider the quadruple 
$(A,B,C,D)$ given by
\beann
A = \bsmat 0 & 1 & 4  \\[1mm] 0 & -1 & 0  \\[1mm] -8 & -6 & 7  \esmat\!\!, \quad
B = \bsmat  0&     2\\[1mm]
     -10  &   0 \\[1mm]
  -2 &   -3 \esmat\!\!, \quad
C = \bsmat 0 & 0 & 1\\[1mm] 0 & 0 & 0 \esmat\!\!, \quad
D = \bsmat -8 & 0 \\[1mm] -2 & 0 \esmat\!\!, \quad
\eeann
which is left and right invertible (so that $\gR^\star=\{0\}$) and has two invariant zeros $z_1=-1$ and $z_2=-16/3$. The null-spaces of $P_{\scriptscriptstyle \Sigma}(-1)$ and $P_{\scriptscriptstyle \Sigma}(-\frac{16}{3})$ are given by
\beann
\ker P_{\scriptscriptstyle \Sigma}(-1)=\ima \bsmat -9\\[1mm]  13 \\[1mm]  0 \\[1mm] \hline 0 \\[1mm] -2  \esmat, \qquad 
\ker P_{\scriptscriptstyle \Sigma}(-\frac{16}{3})=\ima \bsmat -3\\[1mm]  0 \\[1mm]  0 \\[1mm] \hline 0\\[1mm] 8  \esmat.
\eeann
Since $\gR^\star=\{0\}$, we can then select $\hat{V}_g=\bsmat -9 && -3 \\[1mm] 13 && 0 \\[1mm] 0 && 0 \esmat$ and $\hat{W}_g=\bsmat 0 && 0 \\[1mm] -2 && 8 \esmat$. Again, let $\lambda_1=\lambda_2=-2$. We compute basis matrices for $\ker P_{\scriptscriptstyle \Sigma_1}(\lambda_1)$ and $\ker P_{\scriptscriptstyle \Sigma_2}(\lambda_2)$, which are respectively given by
\beann
\bmat{c} V_1 \\ W_1 \emat=\bsmat -3 \\[1mm] 0 \\[1mm] -1 \\[1mm] \hline 0 \\[1mm] 5 \esmat   \quad \text{and} \quad  \bmat{c} V_2 \\ W_2 \emat=\bsmat 73 \\[1mm] 50 \\[1mm] 40 \\[1mm]  \hline 5 \\[1mm] 178 \esmat.
\eeann
Thus, $\gR_1^\star(\lambda_1)=\spanR \{\bmat{cccc} 3 & 0 & 1 \emat^\top \}$ and $\gR_2^\star(\lambda_2)=\spanR \{\bmat{cccc} 73& 50&40 \emat^\top \}$. Condition (\ref{eq:condLhglobgen}) is fulfilled since $\dim(\gV^\star_g+\gR_1^\star(\lambda_1)+\gR_2^\star(\lambda_2))=3$. Taking $H=I_2$, $k_1=2$ and $k_2=0$ (and $K$ is the empty matrix since $\gR^\star=\{0\}$) gives
\beann
F_{\scriptscriptstyle K,H,k_1,k_2,\psi_1} \ns&\ns = \ns&\ns [ \begin{array}{ccc} \hat{W}_g \,H & W_1\,k_1  \end{array} ]\, [ \begin{array}{ccc} \hat{V}_g \,H & V_1\,k_1  \end{array} ]^{-1} \\
\ns&\ns = \ns&\ns \bsmat 0 && 0 && 0  \\[1mm] -2 && 8 && 10 \esmat \bsmat -9 && -3 && -6 \\[1mm] 13 && 0 && 0  \\[1mm] 0 && 0 && -2  \esmat^{-1}= \bsmat 0 && 0 &&  0 \\[1mm] -8/3 && -2 && 3 \esmat,
\eeann
where $\psi_1$ is associated to the choice $k_2=0$.
With this feedback matrix we find $\sigma(A+B\,F_{\scriptscriptstyle K,H,k_1,k_2,\psi_1})=\{-16/3,-1,-2\}$, and the second output component is identically zero.
Notice that if we replace $H$ with any $2 \times 2$ non-singular matrix leads to the same feedback matrix, whereas choosing $k_1=0$ and $k_2=2$ gives
\beann
F_{\scriptscriptstyle K,H,k_1,k_2,\psi_2} \ns&\ns = \ns&\ns [ \begin{array}{ccc} \hat{W}_g \,H &  W_2\,k_2 \end{array} ]\, [ \begin{array}{ccc} \hat{V}_g \,H & V_2\,k_2 \end{array} ]^{-1} \\
\ns&\ns = \ns&\ns \bsmat 0 && 0 && 10 \\[1mm] -2 && 8 && -356 \esmat \bsmat -9 && -3 && 146 \\[1mm] 13 && 0  && 100 \\[1mm] 0 && 0  && 80 \esmat^{-1}= \bsmat 0 && 0 &&  1/8 \\[1mm] -8/3 && -2 && 35/12 \esmat,
\eeann
where $\psi_2$ is associated to the choice $k_1=0$.
With this feedback matrix we still have $\sigma(A+B\,F_{\scriptscriptstyle K,H,k_1,k_2,\psi_2})=\{-16/3,-1,-2\}$, but this time it is the first output component to be equal to zero.
Another solution to Problem~\ref{prob:mono2} is the one in which we insist on forcing both the first and the second components of the output to behave as a single exponential $e^{-2\,t}$, by for example eliminating the second column of $\hat{V}_g$ and $\hat{W}_g$, i.e, 
\beann
F_{\scriptscriptstyle K,H,k_1,k_2} \ns&\ns = \ns&\ns [ \begin{array}{ccc} \hat{W}_g \,H & W_1\,k_1 & W_2\,k_2 \end{array} ]\, [ \begin{array}{ccc} \hat{V}_g \,H & V_1\,k_1 & V_2\,k_2 \end{array} ]^{-1} \\
\ns&\ns = \ns&\ns \bsmat 0 && 0 && 10 \\[1mm] -2 && 10  && -356 \esmat \bsmat -9 && -6 && 146 \\[1mm] 13 && 0 && 100 \\[1mm] 0 && -2 && 80 \esmat^{-1}
\eeann
where $H =\bsmat 1 \\[1mm] 0 \esmat$, and $k_1=k_2=2$. In this case $\sigma(A+B\,F_{\scriptscriptstyle K,H,k_1,k_2})=\{-1,-2\}$, where $-2$ is double. Using $H =\bsmat 0 \\[1mm] 1 \esmat$, and $k_1=k_2=2$ yields
\beann
F_{\scriptscriptstyle K,H,k_1,k_2} \ns&\ns = \ns&\ns [ \begin{array}{ccc} \hat{W}_g \,H & W_1\,k_1 & W_2\,k_2 \end{array} ]\, [ \begin{array}{ccc} \hat{V}_g \,H & V_1\,k_1 & V_2\,k_2 \end{array} ]^{-1} \\
\ns&\ns = \ns&\ns \bsmat 0 && 0 && 10 \\[1mm] 8 && 10  && -356 \esmat \bsmat -3 && -6 && 146 \\[1mm] 0 && 0 && 100 \\[1mm] 0 && -2 && 80 \esmat^{-1}= \bsmat 0 && 1/10 &&  0 \\[1mm] -8/3 && -31/15 && 3 \esmat,
\eeann
and $\sigma(A+B\,F_{\scriptscriptstyle K,H,k_1,k_2})=\{-16/3,-2\}$, where again $-2$ is double. 
\end{example}

\section{Solution to Problem~\ref{prob:mono0}}

In this section, the role played by the ``visible'' eigenvalues $\lambda_1, \lambda_2, \cdots, \lambda_p$ in the solutions to Problem ~\ref{prob:mono0} is investigated. 

\begin{theorem}
\label{th:Lchoice}
Let Assumptions \ref{Ass1} and \ref{Ass2} hold. Let $\dim \gV^\star_g=n-p$. 
Problem~\ref{prob:mono0} is solvable if and only if 

\begin{equation}
\forall\,{S}\in2^{\{1,\ldots,p\}},\;\;\dim\left(\gV_{g}^{\star}+\sum_{j\in{S}}\gR_{j}^{\star}\right)\ge n-p+{\rm {card}}({S}). \label{eq:cond}
\end{equation}
\end{theorem}

\proof
Suppose that (\ref{eq:cond}) is not satisfied. This means that there exists ${S} \in 2^{\{1,\ldots,p\}}$ such that $\dim\left(\gV_{g}^{\star}+\sum_{j\in{S}}\gR_{j}^{\star}\right)<n-p+{\rm {card}({S})}$, which gives $\dim\left(\gV_{g}^{\star}+\sum_{j\in{S}}\gR_{j}^{\star}(\lambda_{j})\right)<n-p+{\rm {card}({S})}$ for any $\lambda_{1},\cdots,\lambda_{p}\in\Lambda_g$,  since by \eqref{eq:gRdec} there holds $\gR_{j}^{\star}(\lambda_{j})\subseteq\gR_{j}^{\star}$ for all $j \in \{1,\cdots, p\}$ and $\lambda_{j} \in \mathbb{R}\backslash \gZ$. In view of Theorem~\ref{th:kchoice}, Problem~\ref{prob:mono2} is never solvable, which implies that Problem~\ref{prob:mono0} does not admit solution.
Sufficiency follows directly from   Corollary \ref{cornew}.
\endproof

\begin{example}
Consider the system in Example \ref{exe0}. If we denote by $\{e_1,\ldots,e_5\}$ the canonical basis in $\gX=\real^5$, we get $\gR^\star_1=\gR^\star_3= \spanR\{e_2,e_3,e_4,e_5\}$, $\gR^\star_2= \spanR\{e_3,e_4,e_5\}$. We recall that $\gV^\star_g=\ima \bsmat -2 && 2/3 && -41/22 && 0 && -1/11 \\[1mm] 0 && 0 && 0 && 1 && 0 \esmat^\top$. In this case, (\ref{eq:cond}) can be written as
\beann
&&\hspace{-4mm}\dim (\gV^\star_g+\gR^\star_j)  \ge  n-p+1 \quad \,j \in \{1,2,3\} \\
&&\hspace{-4mm}\dim (\gV^\star_g+\gR^\star_i+\gR^\star_j)  \ge  n-p+2 \quad \,i,j \in \{1,2,3\}\,\, \text{and $i \neq j$} \\
&&\hspace{-4mm}\dim (\gV^\star_g+\gR^\star_1+\gR^\star_2+\gR^\star_3)  \ge  n-p+3.
\eeann
In the present case, these conditions are satisfied. Indeed, we find $\gV^\star_g+\gR^\star_1=\gV^\star_g+\gR^\star_3=\gX$, the dimension of $\gV^\star_g+\gR^\star_2$ is $4$, and $\gV^\star_g+\gR^\star_2+\gR^\star_3=\gX$.\endex
\end{example}

We now consider the computation of the feedback matrix.

\begin{theorem}
\label{theF2} 
Let Assumptions \ref{Ass1} and \ref{Ass2} hold. Let $\dim \gV^\star_g=n-p$. 
Let the condition in Theorem \ref{th:Lchoice} hold. 
Let $\hat{V}_g$ and $\hat{W}_g$ be constructed as in Lemma \ref{cor2}.
Let $\bsmat V_i(\lambda_i) \\[1mm] W_i(\lambda_i) \esmat$ denote a polynomial basis matrix of least degree for the kernel of $P_{\scriptscriptstyle \Sigma_i}(\lambda_i)=\bsmat A-\lambda_i\,I_n && B \\[1mm] C_{(i)} && D_{(i)} \esmat$ for all $i \in \{1,\ldots,p\}$, where each
$V_i(\lambda_i)$ and $W_i(\lambda_i)$ have $n$ and $m$ rows, respectively. Finally, let $\nu_i$ denote the number of columns of $V_i(\lambda_i)$.\footnote{Notice that $\nu_i$ does not depend on $\lambda_i$ if $\lambda_i \notin \gZ$.} Let $k_i \in \real^{\nu_i}$ denote a parameter vector for all $i \in \{1,\ldots,p\}$. Let
\beann
V_{\scriptscriptstyle K,H,k_1,\ldots,k_p}(\lambda_1, \ldots,\lambda_p) \ns&\ns \!\!\defi  \!\! \ns&\ns [\begin{array}{ccccc} \!\!\!\! \hat{V}_g\diag\{K,H\} \!\!\!\!  &  \!\!\!\! V_1(\lambda_1)\,k_1  \!\!\!\! & \!\!\!\!  \ldots  \!\!\!\! & \!\!\!\!  V_p(\lambda_p)\,k_p  \!\!\!\! \end{array} ] \\
W_{\scriptscriptstyle K,H,k_1,\ldots,k_p}(\lambda_1, \ldots,\lambda_p) \ns&\ns \!\!\defi  \!\! \ns&\ns [\begin{array}{ccccc} \!\!\!\! \hat{W}_g \diag\{K,H\} \!\!\!\!  &  \!\!\!\! W_1(\lambda_1)\,k_1  \!\!\!\! & \!\!\!\!  \ldots  \!\!\!\! & \!\!\!\!  W_p(\lambda_p)\,k_p  \!\!\!\! \end{array} ]
\eeann
Then: \\
{\bf (i)}  the rank of $V_{\scriptscriptstyle K,H,k_1,\ldots,k_p}(\lambda_1, \ldots,\lambda_p)$ is equal to $n$ for almost all $\lambda_i \in \real\setminus \gZ$, for almost all diagonal $K$ and $H$ constructed as in Lemma \ref{cor2} and for  all $k_i \in \real^{\nu_i} \setminus \{0\}$ ($i \in \{1,\ldots,p\}$);\\
{\bf (ii)} The feedback matrices 
\[
F_{\scriptscriptstyle K,H}(\lambda_1, \ldots,\lambda_p)=W_{\scriptscriptstyle K,H,k_1,\ldots,k_p}(\lambda_1, \ldots,\lambda_p)\,V_{\scriptscriptstyle K,H,k_1,\ldots,k_p}^{-1}(\lambda_1, \ldots,\lambda_p)
\]
obtained with $k_i \in \real^{\nu_i}$ and $\lambda_i \in \real_- \setminus \gZ$ ($i \in \{1,\ldots,p\}$) such that $\rank V_{\scriptscriptstyle K,H,k_1,\ldots,k_p}(\lambda_1, \ldots,\lambda_p)=n$ are a solution to Problem~\ref{prob:mono0}.
\end{theorem}
\proof Following essentially the same steps of the proof of Theorem \ref{theF1}, where
now $V_i$, $W_i$ $v_{i,j}$, $w_{i,j}$ are polynomials in $\lambda_i$,  there exist coefficients
$\alpha_{i,j}\in \real$ such that 
\beann
&&\hspace{-0.3cm}\rank [ \begin{array}{c|c|c|c}
\hat{V}_g\,\diag\{K,H\} & \alpha_{1,1}\, v_{1,1}(\lambda_1)+ \ldots+ \alpha_{1,\nu_1}\,  v_{1,\nu_1}(\lambda_1) \end{array} \\
&& \qquad \begin{array}{cccc} \ldots & \alpha_{p,1}\, v_{p,1}(\lambda_p)+ \ldots+ \alpha_{p,\nu_p}\,  v_{p,\nu_p}(\lambda_p)\end{array}] =n.
\eeann
Moreover, the normal rank of $[ \begin{array}{c|c|c|c}
V_g & V_1(\lambda_1) & \ldots & V_p(\lambda_p) \end{array} ]$ is equal to $n$. Since 
\beann
&&\hspace{-0.7cm} [ \begin{array}{c|c|c|c}
\hat{V}_g\,\diag\{K,H\} & \alpha_{1,1}\, v_{i,1}(\lambda_1)+ \ldots+ \alpha_{1,\nu_1}\,  v_{1,\nu_1}(\lambda_1) \end{array} \\
&& \qquad \begin{array}{cccc} \ldots & \alpha_{p,1}\, v_{p,1}(\lambda_p)+ \ldots+ \alpha_{p,\nu_1}\,  v_{p,\nu_p}(\lambda_p)\end{array}]\\
&&\!\!\!\!\!\!\!\!\!\!=[ \begin{array}{c|c|c|c}
 \!\!\! \hat{V}_g \diag\{K,H\} \!\!&\!\!\! V_1(\lambda_1) \!\!\!&\!\! \!\ldots\!\!\! &\!\! V_p(\lambda_p) \!\!\!\!  \end{array} ] \diag\{I_{n-p},\alpha_1,\ldots,\alpha_p\}
\eeann
where $\alpha_i \defi [ \begin{array}{cccc}  \alpha_{i,1} & \ldots & \alpha_{i,p} \end{array} ]^\top$,
the set of parameters $k_1=\alpha_1$, $k_2=\alpha_2$, $\ldots$, $k_p=\alpha_p$ and the set of closed-loop eigenvalues $\lambda_1,\ldots,\lambda_p$ for which 
\beann
[ \begin{array}{c|c|c|c}
 \!\!\! \hat{V}_g\,\diag\{K,H\} \!\!&\!\! V_1(\lambda_1) \!\!&\!\! \ldots\!\! &\!\! V_p(\lambda_p) \!\!\!  \end{array} ]  \diag\{I_{n-p},\alpha_1,\ldots,\alpha_p\}\eeann
loses rank has to satisfy a finite set of linear equations in $k_1,\ldots,k_p$ and $\lambda_1,\ldots,\lambda_p$. 
\endproof

\begin{example}
Consider the quadruple $(A,B,C,D)$ in Example \ref{esprimo}.
For this example, we have
\[
\gR^\star_1=
\ima \bmat{cc}    
\!\!\!\!   1   \!\! & \!\!   0   \!\!\!\! \\     
\!\!\!\!      0  \!\! & \!\!   0   \!\!\!\! \\    
 \!\!\!\!         0   \!\! & \!\!   0   \!\!\!\! \\   
 \!\!\!\!              0   \!\! & \!\!   1    \!\!\!\!  
                  \emat  ,    \quad \gR^\star_2=
\ima \bmat{ccc}   
 \!\!\!\!   0   \!\! & \!\!   1  \!\! & \!\!   0    \!\!\!\! \\    
 \!\!\!\!       1   \!\! & \!\!   0   \!\! & \!\!   0   \!\!\!\! \\ 
 \!\!\!\!       0   \!\! & \!\!   1   \!\! & \!\!   0   \!\!\!\! \\ 
  \!\!\!\!      0   \!\! & \!\!   0   \!\! & \!\!   1  \!\!\!\!       \emat .
\]
Thus, the conditions of Theorem \ref{th:Lchoice} are satisfied, since $\dim (\gV^\star_g+\gR^\star_1)=\dim (\gV^\star_g+\gR^\star_2)=3$ and $\dim (\gV^\star_g+\gR^\star_1+\gR^\star_2)=4$. Polynomial basis matrices of least degree for $\ker P_{\scriptscriptstyle \Sigma_1}(\lambda_1)$ and $\ker P_{\scriptscriptstyle  \Sigma_2}(\lambda_2)$ are given respectively by
\beann
\bmat{c} \!\!\!V_1(\lambda_1) \!\!\! \\ \!\!\! W_1(\lambda_1)\!\!\!  \emat = 
\bsmat 2 & -4 \\[1mm] 0 & 0 \\[1mm]  0 & 0 \\[1mm] 0 & 4(\lambda_1-1)  \\[1mm] \hline \\ 1 & -2  \\[1mm]  \lambda_1-1 & 0  \\[1mm] 0 & (\lambda_1-1)^2 \esmat\!\!\!, 
\eeann
and 
\beann
\bmat{c} \!\!\! V_2(\lambda_2) \!\!\! \\ \!\!\! W_2(\lambda_2) \!\!\! \emat =\bsmat 
2(\lambda_2+2) & -4(\lambda_2+2)  \\[1mm] 
-4 & 8 \\[1mm]  
2(\lambda_2+2) & -4(\lambda_2+2)  \\[1mm] 
0 & 4(2\,\lambda_2+3)(\lambda_2+2)  \\[1mm] 
\hline \\
-(\lambda_2+3)(\lambda_2+2) & 2(\lambda_2+3)(\lambda_2+2)  \\[1mm] 
(\lambda_2+2)(2\,\lambda_2+3) & 0 \\[1mm] 
0 & (2\,\lambda_2+3)(\lambda_2-1)(\lambda_2+2)\esmat
\eeann
It is easily seen that considering the parameter vectors $k_1=\bsmat k_{11} \\[1mm] k_{12} \esmat$ and $k_2=\bsmat k_{21} \\[1mm] k_{22} \esmat$ we have
\beann
&& \hspace{-0.7cm} [ \begin{array}{c|c|c}
\hat{V}_g \,H & V_1(\lambda_1)\,k_1 & V_2(\lambda_2)\,k_2 \end{array} ] = \\
&& \bmat{cccc}  1 & 0 & 2\,k_{11}-4\,k_{12} & 2\,k_{21}(\lambda_2+2)-4\,k_{22}(\lambda_2+2) \\
2 & 0 & 0 & -4\,k_{21}+8\,k_{22} \\ 1 & 0 & 0 & 2(\lambda_2+2)\,k_{21}-4\,k_{22} (\lambda_2+2) \\
0 & 1  & 4\,k_{12}(\lambda_1-1) & 4\,k_{22}(2\,\lambda_2+3)(\lambda_2+2) \emat.
\eeann
The determinant of this matrix is the polynomial $8\,(2\,k_{12}-k_{11})(k_{21}-2\,k_{22})(3+\lambda_2)$, which is different from zero for almost all parameters $k_{ij}$ and eigenvalues $\lambda_1,\lambda_2$.
Since $H$ is full column-rank, the feedback matrices do not depend on $H$. Thus, the set of all the feedback matrices that solve Problem~\ref{prob:mono0} are given in parameterized form as
\beann
&&\hspace{-0.9cm}F(\lambda_1,\lambda_2)  =\bsmat
0 && 0 && k_{11}-2\,k_{12} && (\lambda_2+3)(\lambda_2+2)(2\,k_{22}-k_{12}) \\[1mm]
-\frac{3}{2} && \frac{1}{2} && k_{11}(\lambda_1-1) && k_{12}\,(\lambda_2+2)(2\lambda_2+3)  \\[1mm]
0 && -1 && k_{12}\,(\lambda_1-1)^{2} && k_{22}\,(2\,\lambda_2+3)(\lambda_2-1)(\lambda_2+2)\esmat \\
&& \qquad \times \bsmat 1 && 0 && 2\,k_{11}-4\,k_{12} && 2\,k_{21}(\lambda_2+2)-4\,k_{22}(\lambda_2+2)  \\[1mm]
2 && 0 && 0 && -4\,k_{21}+8\,k_{22} \\ 1 && 0 && 0 && 2(\lambda_2+2)\,k_{21}-4\,k_{22} (\lambda_2+2)  \\[1mm]
0 && 1  && 4\,k_{12}(\lambda_1-1) && 4\,k_{22}(2\,\lambda_2+3)(\lambda_2+2) \esmat^{-1},
\eeann
where $k_{ij}$ and $\lambda_2$ are such that the determinant $8\,(2\,k_{12}-k_{11})(k_{21}-2\,k_{22})(3+\lambda_2)$ is not zero.\endex
\end{example}

The following result generalises the conditions obtained above for Problem~\ref{prob:mono0} to the case in which $\dim \gV^\star_g\ge n-p$.

\begin{theorem}
\label{the new}
 Let Assumptions \ref{Ass1} and \ref{Ass2} hold. Let $h \defi \dim \gV^\star_g$.
Problem~\ref{prob:mono0} is solvable if and only if 
\begin{equation} 
\dim\left(\gV_{g}^{\star}+\sum_{j\in{S}}\gR_{j}^{\star}\right)\ge n-p+{\rm {card}({S})}
\end{equation}
holds true for all ${S}\in \{\mathfrak{S} \in 2^{\{1,\ldots,p\}}\,|\;\; \mathrm{card}\,{\mathfrak{S}}>h-(n-p)\}$.
\end{theorem}

Observe also that the necessary and sufficient condition in Theorem \ref{the new}  can be alternatively written as 
\beann
\dim(\gV^\star_g+\gR^\star_{\nu_1}+\ldots+\gR^\star_{\nu_l})\ge (n-p)+l 
\eeann
for every $1 \le \nu_1 <  \ldots < \nu_l \le p$ and $l \in \{h-(n-p),\ldots,p\}$.

\begin{remark}
While we have shown an example of a non-minimum phase systems in which the problem of obtaining a monotonic response from any initial conditions can be solved, it is well known that a SISO strictly proper system with real non-minimum phase zeros cannot be monotonic as undershoot must occur, \cite{Middleton-91}. This fact also follows as a particular 
case of Theorem \ref{the new}. Indeed, the condition $\dim \gV_{g}^{\star}\ge n-p$ follows from (\ref{eq:cond}) when $S=\emptyset$. Such condition is never satisfied for SISO strictly proper non-minimum phase systems.  In fact, we have the inequalities $\dim  \gV_{g}^{\star} \le \dim \gV^\star \le \dim (\ker C) = n-1$, where the first can be an equality only if the system is minimum-phase. 
\end{remark}

\begin{example}
Consider the quadruple given in Example \ref{esaggiunto2}.
Here we have
$\gR^\star_1=
\ima \bsmat 1 & 0 \\[1mm] 0 & 0 \\[1mm] 0 & 1 \esmat$ and     $\gR^\star_2=\{0\}$, so that the condition of Theorem \ref{th:Lchoice} is not satisfied. However, since $\dim \gV^\star_g=2>1=n-p$, the solvability condition is the one given in Theorem \ref{the new}, which in this case reduces to $\dim(\gV^\star_g+\gR^\star_1+\gR^\star_2)=n$. This condition is satisfied, so that the problem is solvable. 
A polynomial basis matrix of least degree for $\ker P_{\scriptscriptstyle \Sigma_1}(\lambda)$ is given by
\beann
\bsmat V(\lambda) \\[1mm]  W(\lambda) \esmat = 
[\begin{array}{ccc|cc}  \!  \! -16 \!  &  \! 0  \!  & \!   -2\,(\lambda+10) \! &  \! \frac{\lambda^2-108}{4} \! & \! 3\,(\lambda+22) \!  \!  \end{array}]^\top\!\!\!.
\eeann
One can directly check that with e.g. $H=I_2$ the rank of 
\beann
[ \begin{array}{c|c} V_g\,H & V(\lambda)\,k \end{array} ]= \bmat{cc|c}
8 & 0 & -16\,k \\ 0 & 1 & 0 \\ 3 & 0 & -2\,k\,(\lambda+10) \emat
\eeann
is equal to $3$ for any $k \neq 0$ and any $\lambda \neq -7$, and the  feedback matrices that solve the problem are parameterised in $\lambda$ as
\beann
F_{\scriptscriptstyle K,H,k}(\lambda)=\bmat{ccc}   
\!\!\!\!  \frac{59}{8}   \!\! & \!\!    \frac{3}{8} \!\! & \!\! \frac{\lambda^2-108}{4}\,k  \!\!\!\! \\        
\!\!\!\!  -\frac{45}{2}    \!\! & \!\!   \frac{5}{4}  \!\! & \!\!3\,k\,(\lambda+22)    \!\!\!\! 
      \emat\!\!\!\bmat{ccc}
8 & 0 & -16\,k \\ 0 & 1 & 0 \\ 3 & 0 &\!\!\! -2\,k\,(\lambda+10) \!\!\!\emat^{-1}
\eeann
for $k \neq 0$ and $\lambda \neq -7$.
For example, choosing $\lambda=-10$ gives
\beann
F=F_{\scriptscriptstyle K,H,k}(-10)=\bmat{ccc}  \frac{1}{8} & \frac{3}{8} & \frac{17}{8} \\ -\frac{9}{4} & \frac{5}{4} & -\frac{3}{2} \emat.
\eeann
 Notice that $\sigma(A\!+\!B\,F)=\{-7,-9,-10\}$ and $(C\!+\!D\,F)\,\bsmat 8 & 0 & -8 \\[1mm] 0 & 1 & 0\\[1mm] 3 & 0 & 0 \esmat=\bsmat 0 & 0 & \Big|& 24 \\[-3mm] 0 & 0 & \Big| & 0 \esmat$. With this feedback matrix the second output instantaneously tracks the reference.\endex
\end{example}

\section*{Concluding remarks}
In this paper, the problem of achieving a monotonic step response from any initial condition has been addressed for the first time for LTI MIMO systems.
 This new approach opens the door to a range of developments that for the sake of conciseness cannot be addressed in this paper, but that we briefly discuss:
 \begin{itemize}
 \item In the case where global monotonicity cannot be achieved, it is important to find structural conditions ensuring that every component of the tracking error consists of the sum of at most two, three, or more closed-loop modes. In such case, even if the response is not globally monotonic, it is still monotonic starting from suitable initial conditions. An important issue is the characterisation of the regions of the state space where the initial state must belong to guarantee that the response can be made monotonic;
 \item We have parameterised the set of feedback matrices that solve the problem of obtaining a monotonic step response from an arbitrary initial condition. A second relevant problem
 involves the use of the method in \cite{SICON} to the end of exploiting the remaining degrees of freedom in the parameterisation in order to compute the state feedback that achieves a globally monotonic step response and which at the same time delivers a robust closed-loop eigenstructure, by ensuring that the closed-loop eigenvalues are rendered as insensitive to perturbations in the state matrices as possible. 
This task can be accomplished by
 obtaining a feedback matrix that minimises the Frobenius condition number of the matrix of closed-loop eigenvectors, which is a commonly used robustness measure. 
 The problem of obtaining a feedback matrix with minimum gain measure can be handled in a similar way, by minimising 
 the Frobenius norm of the feedback matrix. 
 \end{itemize}

%

\end{document}